\documentclass[12pt]{amsart}

\usepackage{amsmath,amssymb}
\usepackage{amsfonts}
\usepackage{amsthm}
\usepackage{latexsym}
\usepackage{graphicx}
\usepackage{epsfig}

\newtheorem{theorem}{Theorem}[section]
\newtheorem{lemma}[theorem]{Lemma}

\theoremstyle{definition}

\newtheorem{example}{Example}[section]

\theoremstyle{remark}
\newtheorem{claim}{Claim}[section]
\newtheorem{remark}{Remark}[section]

\makeatletter \@addtoreset{equation}{section} \makeatother

\def\ddt{\frac{d}{dt}}

\begin{document}

\title{Cross Curvature Flow on Locally Homogeneous Three-manifolds (II) }

\author{Xiaodong Cao$^*$}
\thanks{$^*$Research
partially supported by the Jeffrey Sean Lehman Fund from Cornell
University}

\address{Department of Mathematics,
  Cornell University, Ithaca, NY 14853}
\email{cao@math.cornell.edu, lsc@math.cornell.edu}

\author{Laurent Saloff-Coste$^\flat$}
\thanks{$^{\flat}$Research
partially supported by NSF grant no. \# DMS 0603886}


\renewcommand{\subjclassname}{%
  \textup{2000} Mathematics Subject Classification}

\date{May 15th,  2008}

\maketitle

\markboth{Xiaodong Cao,  Laurent Saloff-Coste} {Cross Curvature
Flow on Locally Homogeneous Three-manifolds (II)}

\begin{abstract}  In this paper, we study the positive cross
curvature flow on locally homogeneous $3$-manifolds. We describe
the long time behavior of these flows. We combine this with
earlier results concerning the asymptotic behavior of the negative
cross curvature flow to describe the two sided behavior of maximal
solutions of the cross curvature flow on locally homogeneous
$3$-manifolds. We show that, typically, the positive cross
curvature flow on locally homogeneous $3$-manifold produce an
Heisenberg type sub-Riemannian geometry.
\end{abstract}

\section{Introduction}
\subsection{Evolution equations on homogeneous manifolds}
Hamilton's Ricci flow (\cite{H3}) is the best known example of a
geometric evolution equation.  One of the aims  of such flows is
to obtain metrics with special properties (in the case of the
Ricci flow, Einstein metrics). Special cases arise when the metric
is invariant under a group of transformations and this property is
preserved by the flow. In particular, if the group of isometries
of the original Riemannian structure is transitive, then the
geometric evolution equation reduces to an ODE in the tangent
space of an arbitrary fixed origin. The Ricci flow on locally
homogeneous $3$-manifolds was analyzed in \cite{isenberg1} and on
some homogeneous $4$-manifolds in \cite{isenberg2}.

The cross curvature flow, or (XCF), was introduced by Chow and
Hamilton \cite{chowhamiltonxcf} and depends on the choice of a
sign (see Section 1.3 below) leading to two flows: (+XCF) and
(-XCF). Chow and Hamilton conjectured that for any compact
$3$-manifold that admits a metric with negative sectional
curvature, the normalized positive cross curvature flow, started
at such a metric, exists for all time and converges to a
hyperbolic metric. In \cite{cnsc1}, we study the asymptotic
behavior of the negative cross curvature flow (-XCF) on
homogeneous $3$-manifolds. In this paper, we complement the
results of \cite{cnsc1} by studying the asymptotic behavior of the
positive cross curvature flow, or (+XCF), on homogeneous
$3$-manifolds. In the homogeneous case, the local existence is not
an issue and  the negative and positive cross curvature flows can
be seen as the same basic flow (say, (-XCF)) run either in the
forward  or backward direction. Putting together the results of
\cite{cnsc1} and of this paper, we will thus describe both the
forward and backward asymptotic behaviors of the maximal solution
of (-XCF) through any given metric $g_0$ on any locally
homogeneous $3$-manifold.

Although we will give much more precise statements,
the spirit of the main results proved in this paper is captured in the
following theorem.

\begin{theorem}\label{th-main}
Let $(M, g_0)$ be a complete locally homogeneous $3$-manifold
(compact or not). Let $g^b(t), t\in [0,T_b)$ be the maximal
solution of the positive cross curvature flow {\em (+XCF)} with
initial metric $g_0$. Let $d(t)$ denote the distance function on
$(M,g^b(t))$. Assume that $g_0$ is generic among all locally
homogeneous metrics on $M$. Then \begin{itemize} \item either
$T_b=\infty$ and $g^b(t)= e^{\lambda t} g_0$ for some $\lambda\in
\mathbb R$ (i.e., the cross curvature tensor of $g_0$ is equal to
$\lambda g_0$), \item or $T_b<\infty$ and  there exists a function
$r(t):[0,T_b)\to (0,\infty)$ such that the metric spaces
$(M,r(t)d(t))$ converge uniformly to a sub-Riemannian metric space
$(M,d(T))$ whose tangent cone at any  $m\in M$ is the Heisenberg
group $\mathbb H_3$ equipped with its natural sub-Riemannian
metric.
\end{itemize}
\end{theorem}
In this paper, we prove this statement in all cases except for
manifolds covered by $\widetilde{\mbox{SL}(2,\mathbb R)}$ (the
$\mbox{SL}(2,\mathbb R)$ case). In the $\mbox{SL}(2,\mathbb R)$
case, we only prove a slightly weaker result (see Theorem
\ref{th-sl2}). To obtain a proof of Theorem \ref{th-main} in the
 $\mbox{SL}(2,\mathbb R)$ case, one actually needs some additional
 information. This additional information is
 obtained in \cite{cgsc08} by using a different type of argument.
 See the comment at the end
 of the proof of Theorem \ref{th-sl2}.

\begin{remark}
The first case only occurs for homogeneous $3$-manifolds covered
by $\mathbb R\sp 3$, $\mathbb H\sp 3$, $\mathbb S\sp 2 \times
\mathbb R$ and $\mathbb H\sp 2 \times \mathbb R$. Moreover, in
those cases, $g^b(t)= e^{\lambda t} g_0$ for any homogeneous
$g_0$.
\end{remark}
We find it quite striking that the asymptotic behavior of (+XCF),
i.e., the backward behavior of (-XCF), is essentially the same in
all cases, for generic homogeneous metrics. In a companion paper
\cite{cscbrf}, we show that the same ``universal'' behavior holds
for the backward Ricci flow on homogeneous $3$-manifolds. This
contrast with the very different behavior observed in the forward
direction. See \cite{isenberg1,cnsc1} and the various more precise
statements given below.

\subsection{The cross curvature tensor on $3$-manifolds}

On a $3$-dimensional Riemannian manifold $(M,g)$, let $Rc$ be the
Ricci tensor and $R$ be the scalar curvature. The Einstein tensor
is defined by $E=Rc-\frac{1}{2}Rg$. Its local components are
$E_{ij}=R_{ij}-\frac{1}{2}Rg_{ij}.$  Raising the indices, define
$P^{ij}=g^{ik}g^{jl}R_{kl}-\frac{1}{2}Rg^{ij},$ where $g^{ij}$ is the
inverse of $g_{ij}$. Let $V_{ij}$ be the inverse of $P^{ij}$ (if
it exists). The cross curvature tensor is (see \cite{chowhamiltonxcf})
$$h_{ij}=\left(\frac{\det P^{kl}}{\det g^{kl}}\right)V_{ij}.$$
Assume that computations are done in an orthonormal frame where
the Ricci tensor, and thus also the cross curvature tensor, is
diagonal. If the principal sectional curvatures are $k_1,k_2,k_3$
($k_i=K_{jkjk}$, circularly) so that
$R_{ii}=k_j+k_l$, circularly, then
\begin{equation} \label{C=kk}
h_{ii} = k_jk_l.
\end{equation}
Notice that this definition works even when some of the sectional
curvatures vanish.

\subsection{The cross curvature flows}
In \cite{chowhamiltonxcf}, Chow and Hamilton define the cross
curvature flow on $3$-manifolds starting from a metric with
either positive sectional
curvature or negative sectional curvature. More precisely,
if $\epsilon =\pm 1$ is the sectional curvature sign
(assumed to be well defined) of the metric $g_0$,
the cross curvature flow starting from $g_0$ is the
solution of $$\left\{\begin{array}{l}\frac{\partial}{\partial
t}g=-2\epsilon h\\
g(0)=g_0.\end{array}\right.$$ In these circunstances,
the local existence of the flow was
proved in \cite{buckland}.

Locally homogeneous manifolds seldom have sectional curvatures
that are all of the same sign. In dimension 3, positive sectional
curvature is only possible on locally homogeneous manifolds
covered by the sphere $\mbox{SU}(2)$. Negative sectional curvature
occurs only on hyperbolic $3$-manifolds. All other locally
homogeneous closed Riemannian $3$-manifolds are either flat or
have some positive sectional curvature \cite[Theorem
1.6]{milnor76}. Thus the definition above is not really practical
for our purpose as far as a choice of sign is concerned and it is
natural to investigate both the positive and the negative cross
curvature flows defined by
$$\left\{\begin{array}{l}
\frac{\partial}{\partial t}g=2 h\\
g(0)=g_0.\end{array}\right.\leqno{(+\mbox{XCF})}$$ and
$$\left\{\begin{array}{l}
\frac{\partial}{\partial t}g= -2 h\\
g(0)=g_0.\end{array}\right. \leqno{(-\mbox{XCF})}$$ In fact,
starting from a initial metric $g_0$ on a locally homogeneous
$3$-manifold, let $g^f(t)$, $t\in [0,T_{f})$ be the maximal
forward solution of the (-XCF), and $g^b(t)$, $t\in [0,T_b)$ be
the maximal forward solution of the (+XCF). Now, for $t\in
I=(-T_b,T_f)$, set
$$g(t)=\left\{ \begin{array}{cc}
g^f(t)& \mbox{ for } t\in [0,T_f)\\
g^b(-t) & \mbox{ for } t\in (-T_b,0].\end{array}\right.$$ It is
easy to see that, by construction, $g(t), t\in I$, is a maximal
solution of (-XCF) passing through $g_0$ at time $t=0$. One of the
goals of our study is to describe the behavior of these maximal
solutions in the forward and backward directions. The forward
direction is treated in \cite{cnsc1} and the present paper is
devoted to the backward direction which, of course, is the same as
the forward direction for (+XCF).

\subsection{Normalizations}
Let $g(t), t\in I$, be a maximal solution of the (-XCF). By
renormalization of $g(t)$, we mean a family
$\widetilde{g}(\widetilde{t}), \widetilde{t}\in \widetilde{I},$
obtained by a change of scale in space and a change of time, that
is
$$\widetilde{g}(\widetilde{t})=\psi (t) g(t),~
\widetilde{t}=\int_{0}^t \psi^2(s) ds.$$

Set $\widetilde{\psi} (\widetilde{t})=\psi (t)$, then we have
$$\frac{\partial \widetilde{g}}{\partial
\widetilde{t}}=-2\widetilde{h}+\left(\frac{d}{d \widetilde{t}} \ln
\widetilde{\psi}\right) \widetilde{g},$$ where $\widetilde{h}$ is
the cross curvature tensor of $\widetilde{g}$.

On compact manifolds, it is customary to take $\ddt \ln
\psi=\frac{2}{3} \overline{h}$, where $\overline{h}=\frac{\int
tr(h) d\mu}{\int d\mu}$ is the average of the trace of the cross
curvature because the volume of the metric $\widetilde{g}$ is then
constant. In this paper, we will consider some different
normalizations, for instance, keeping the diameter constant.
Moreover, we will not worry about the time change associated above
with a re-scaling by $\psi$. Given a solution $g$ of (-XCF), we
will be interested in finding re-scaling $\bar{g}(t)=\phi(t)g(t)$
such that the asymptotic behavior of the metric space
$(M,\bar{g}(t))$ is described by a model having the largest
possible dimension (i.e., minimum collapse). Although making this
precise could possibly involve some difficulties in general, in
our specific examples, what it means will be quite obvious.

\subsection{Convergence of metric spaces}\label{sec-convd}
We refer the reader to \cite{BBI} for an introduction to and more details on
the notions discussed briefly in this section.
We start with the most basic (and naive) notion of convergence:
the uniform convergence of a family of metric spaces
$(X,d_t)$, $t\in (0,T)$, towards a metric space $(X,d_T)$
when all metric structures are defined on the same
topological space $X$. By definition, this uniform convergence
(which we will encounter frequently below), is simply the
convergence of the functions $d_t$ to $d_T$ on
$X\times X$, uniformly on compact sets, as $t$ tends to $T$.
\begin{example} \label{e-gg}
Let $g_t$, $t\in (0,T]$,  be Riemannian metrics
on a manifold $M$, equipped with an auxiliary Riemannian metric
$g_0$. Let $d_t$ be the corresponding distance functions on $M$.
Assume that there is a Riemannian metric $g_T$ on $M$ such that
for any compact $K\subset M$,
\begin{equation}\label{gg}
\lim_{t\rightarrow T} \max_{x\in K}\max_{u\in T_x:
g_0(u,u)\le 1}| g_t(u,u)-g_T(u,u)|=0.\end{equation}
Then the distance $d_t$
converge uniformly on compact sets to $d_T$.
\end{example}

The Lipschitz distance between two metric spaces is the infimum of
$\ln \mbox{dil}(f)+\ln \mbox{dil}(f^{-1})$ where $f$ is a
bi-lipschitz homeomorphism between $X$ and $Y$ and
\begin{equation}\label{def-dil}
\mbox{dil}(f)=
\sup_{x,x'\in X}\{d_Y(f(x),f(x'))/d_X(x,x')\}.\end{equation}

\begin{example}Let $M$ and $g_t$, $t\in [0,T]$ be as in Example \ref{e-gg}.
If (\ref{gg}) holds and $M$ is compact then the Lipschitz distance between
$(M,d_t)$ and $(M,d_T)$ tends to zero. If $M$ is not compact, this is not
necessarily the case because Lipschitz convergence implies
some global control. However, if one assumes instead  that
\begin{equation}\label{ggg}
\lim_{t\rightarrow T} \max_{x\in M}\max_{u\in T_x:
g_0(u,u)\le 1}| g_t(u,u)-g_T(u,u)|=0\end{equation}
then  the Lipschitz distance between  $(M,d_t)$ and $(M,d_T)$ tends to zero.
\end{example}

Recall that for two subsets $A,B$ of a metric space $Z$,
$$d^Z_H(A,B)=\inf \{\epsilon>0: A\subset B_\epsilon \mbox{ and }
B\subset A_\epsilon\}$$ where $A_\epsilon$ is the
$\epsilon$-neighborhood of $A$ in $Z$. The Hausdorff distance
between two metric spaces $X,Y$ is the infimum of the numbers
$d^Z_H(f(X),g(Y))$ for all metric spaces $Z$ and all isometric
embeddings $f,g$ of $X,Y$ into $Z$. Finally, a sequence of pointed
metric spaces $(X_n,p_n)$ converges in the Gromov-Hausdorff sense
to a pointed metric space $(X,p)$ if, for every $r,\epsilon>0$,
there is a map $f: B_{X_n}(p_n,r)\rightarrow X$ (not necessarily
continuous) such that  $f(p_n)=p$, $\mbox{dis}(f)<\epsilon$,
 $B_X(p,r-\epsilon)\subset [f(B_{X_n}(,p_n,r))]_\epsilon$.
 Here, if $f$ is a map from $X$ to $Y$,
 $\mbox{dis}(f)=\sup_{x,x'\in X}|d_Y(f(x),f(x'))-d_X(x,x')|$.
For length spaces (and we will deal only with length spaces),
this is equivalent to say that the balls $B_{X_n}(p_n,r)$ converge
in the Gromov-Hausdorff sense to $B(p,r)$, for all each $r>0$.

\begin{example} Convergence in the Gromov-Hausdorff sense allows for
dimensional collapse. For instance, the pointed cylinder
$(\mathbb R\times S^1,(0,0), d_t)$,
$d_t$ being the distance associated with
$g_t= (dx)^2+ t^{-1}(d\theta)^2$  converges in
the Gromov-Hausdorff sense as $t$ tends to infinity
to $(\mathbb R, 0)$ equipped with its usual metric.

Note that the metric spaces $(\mathbb R^2,(0,0), d_t)$, $d_t$ associated with
$g_t=(dx)^2+t^{-1}(dy)^2$, are all isometric and thus indistinguishable in
terms of the Lipschitz or Gromov-Hausdorff distances.
 \end{example}

\begin{example}[Tangent cones] Given a pointed metric space $(X,d,p)$,
we call tangent cone at $p$, any pointed metric space $(X_0,d_0,p_0)$
which appears as a Gromov-Hausdorff limit of (a subsequence of)
the family of pointed metric spaces $(X,td,p)$, $t$ tending to $0$.
For any pointed Riemannian $n$-manifold $(M,g,p)$, the tangent cone at $p$
exists, is unique, and equals the Euclidean $n$-space.
\end{example}

\subsection{Locally homogeneous  $3$-manifolds}
By classical arguments, the study of the Ricci or cross curvature
flow on a locally homogeneous manifold reduces essentially to the
study of the same flow on the universal cover. In dimension $3$
there are 9 possibilities for the universal cover, four of which
are essentially trivial as far as the cross curvature flow is
concerned. These four easy cases are : $\mathbb R^3$ (flat
metrics), $\mathbb H^3$ (hyperbolic metric), $\mathbb S^2\times
\mathbb R$ and $\mathbb H^2\times \mathbb R$. See \cite{cnsc1}. In
the remaining 5 cases the universal cover is itself a group that
act transitively on the manifold. This paper focusses on these
five cases which are:
 $\mbox{SU}(2)$, $\widetilde{\mbox{SL}(2,\mathbb R)}$; Heisenberg;
$E(1,1)=\mbox{Sol}$ (the group of isometry of plane with a flat
Lorentz metric); $\widetilde{E(2)}$ (the universal cover of the
group of isometries of the Euclidean plane). See \cite{isenberg1}
or \cite{cnsc1} for a more detailed discussion.

Assume that $\mathfrak g$ is a $3$-dimensional real Lie unimodular
algebra equipped with an oriented Euclidean structure. According
to J. Milnor, \cite{milnor76},
 there exists a (positively oriented) orthonormal basis
$(e_1,e_2,e_3)$ and reals $\lambda_1,\lambda_2,\lambda_3$ such
that the bracket operation of the Lie algebra has the form
$$[e_i,e_j]=\lambda_k e_k  \;\;\;\mbox{(circularly in $i,j,k$)}.$$
Milnor shows that such a basis diagonalizes the Ricci tensor and
thus also the cross curvature tensor. If $f_i= a_ja_k e_i$ with
nonzero $a_i,a_j,a_k\in \mathbb R$, then $[f_i,f_j]= \lambda_k
a_k^2 f_k$ (circularly in $i,j,k$). Using the choice of
orientation, we may assume that at most one of the $\lambda_i$ is
negative and then, the Lie algebra structure is entirely
determined by the signs (in $\{-1,0,+1\}$) of
$\lambda_1,\lambda_2,\lambda_3$ as follows:
$$\begin{array}{cccl}
+&+&+&  \;\mbox{SU}(2)\\
+&+&-& \; \mbox{SL}(2,\mathbb R) \\
+&+&0& \;  E(2) \; \mbox{(Euclidean
motions in $2D$)} \\
+&0&-&  \;E(1,1) \; \mbox{(also called Sol)}\\
+&0&0& \;  \mbox{Heisenberg group} \\
0&0&0& \;  \mathbb R^3
\end{array}$$
In each case, let $\epsilon=(\epsilon_1,\epsilon_2,\epsilon_3)\in
\{-1,0,+1\}^3$ be the corresponding choice of signs. Then, given
$\epsilon$ and an Euclidean  metric $g_0$ on the corresponding Lie
algebra, we can choose a basis $f_1,f_2,f_3$ (with $f_i$ collinear
to $e_i$ above) such that
\begin{equation}\label{MF}
[f_i,f_j]= 2\epsilon_k f_k \;\;\;\mbox{(circularly in $i,j,k$)}.
\end{equation}
We call $(f_i)_1^3$ a Milnor frame for $g_0$. The metric, the
Ricci tensor and the cross curvature tensor are all diagonalized
in this basis and this property is obviously maintained throughout
either the Ricci flow or cross curvature flow. If we let
$(f^i)_1^3$ be the dual frame of $(f_i)_1^3$, the metric $g_0$ has
the form
\begin{equation}\label{g_0}
g_0= A_0 f^1\otimes f^1 +B_0f^2\otimes f^2+ C_0f^3\otimes
f^3.\end{equation} Assuming existence of the flow $g(t)$ starting
from $g_0$, under either the Ricci flow or the cross curvature
flow (positive or negative), the original frame $(f_i)_1^3$ stays
a Milnor frame for $g(t)$ along the flow and $g(t)$ has the form
\begin{equation}\label{gflow}
g(t)= A(t) f^1\otimes f^1 +B(t) f^2\otimes f^2+ C(t) f^3\otimes
f^3.\end{equation} It follows that these flows reduce to ODEs in
$(A,B,C)$. Given a flow, the explicit form of the ODE depends on
the underlying Lie algebra structure. With the help of the
curvature computations done by Milnor in \cite{milnor76}, one can
find the explicit form of the equations for each Lie algebra
structure. The Ricci flow case was treated in \cite{isenberg1}.
The computations of the ODEs corresponding to the cross curvature
flow are presented in \cite{cnsc1} and will be used below to study
the asymptotic behavior of (+XCF).

\subsection{Sub-Riemannian Geometry}
The notion of sub-Riemannian geometry can be described from
several equivalent but different viewpoints. The simplest is
perhaps to start with a family of smooth vector fields $\mathcal
X=\{X_i,i\in \{1,\dots,k\}\}$ on a manifold $M$ with the property
that these fields, together with their brackets of all orders,
span the tangent space at each point of $M$. This is often called
H\"ormander's condition (H\"ormander proved that the associated
sum of squares $L=\sum X_i^2$ is hypoellitic). Given a family $\mathcal X$,
one
defines a distance on $M$ as follows. A vector $u$ in the tangent
space $T_x$ at $x$ is subunit  if $u=\sum_1^k a_iX_i(x)$ with
$\sum |a_i|^2\le 1$. An absolutely continuous curve $\gamma: [0,T]
\rightarrow M$ is subunit if $\dot\gamma (t) $ is subunit for each
$t\in [0,T]$ (in particular, this means that $\dot\gamma(t)$
belongs to the span of $\mathcal X$ at $\gamma(t)$, i.e., is
horizontal). The distance $d_\mathcal X(x,y)$ is the infimum of
$T$ such that there exists an absolutely continuous subunit curve
$\gamma: [0,T]\rightarrow M$ with $\gamma(0)=x$ and $\gamma(T)=y$.

Another equivalent definition starts with a distribution $H$, that
is to say, a sub-bundle of the tangent bundle, together with a
fiber inner product on this sub-bundle. This easily leads to the
notion of length of a horizontal curve (i.e., a curve that stays
tangent to the sub-bundle). In this case, the
H\"ormander condition is expressed using a local frame for $H$.

A third equivalent definition is based on the choice of a
symmetric non-negative $(0,2)$-tensor $Q$ (a possibly degenerate
inner product on the co-tangent bundle). This
defines  a sub-bundle of the
tangent bundle by setting
$$H_x=\{u\in T_x: \sup_{\alpha: Q(\alpha,\alpha)\le
1}\alpha(u)<\infty\}$$ and induces an inner product on $H_x$ in
the obvious way. Again, H\"ormander's condition can be expressed
using a local frame for $H$. The link between the first
presentation and the third is simply that, given a H\"ormander
family $\mathcal X$,
 \begin{equation}\label{def-Q}
Q(\alpha,\alpha)=\sum_1^k |\alpha(X_i)|^2.\end{equation}

The most basic result of sub-Riemannian geometry is that, assuming
H\"ormander's condition, the associated sub-Riemannian distance
defines the original topology of the manifold $M$. See Chow's
theorem in \cite[Ch. 2]{Mont}. More generally, we refer the reader
to \cite{Mont} for a detailed introduction to sub-Riemannian
geometry.

Most relevant to the present paper is the fact that sub-Riemannian
metrics can easily appear as limit of Riemannian metrics as
explain in the following example.

\begin{example} Let $M$ be a manifold, equipped
with a family of Riemannian metric $g_t$, $t\in [0,T)$ (we will
use $g_0$ as a reference metric here). For each $t\in [0,T)$, we
let $d_t$ be the corresponding distance function. Each $g_t$
induces a symmetric positive definite $(0,2)$-tensor $Q_t$. Now,
the existence of a Riemannian metric  $g_T$ such that
$$\lim_{t\rightarrow T} \max_{x\in K}\max_{u\in T_x:\atop
g_0(u,u)\le 1}| g_t(u,u)-g_T(u,u)|=0$$
is obviously equivalent to the the existence of a symmetric positive
definite $(0,2)$-tensor $Q_T$ such that
$$\lim_{t\rightarrow T} \max_{x\in K}\max_{u\in T'_x:\atop
Q_0(u,u)\le 1}| Q_t(u,u)-Q_T(u,u)|=0.$$
However, in general, it is well possible that there exists a
symmetric non-negative $(0,2)$-tensor $Q_T$ such that
$$\lim_{t\rightarrow T} \max_{x\in K}\max_{u\in T'_x:\atop
Q_0(u,u)\le 1}| Q_t(u,u)-Q_T(u,u)|=0$$
even if the metrics $g_t$ do not have a well defined  finite limit.
In that case, if the limiting $(0,2)$-tensor $Q_T$ turns out to satisfy
H\"ormander's condition then the metric space $(M,d_t)$ converges
uniformly on compact sets to the sub-Riemannian metric space $(M,d_T)$
where $d_T$ is the sub-Riemannian distance function associated with $Q_T$.
See, e.g., \cite{J-SC}.
\end{example}

The case of left-invariant sub-Riemannian structures on Lie
groups is somewhat simpler than the general case and extremely
natural. Recall that the Lie algebra $\mathfrak g$ of a connected
Lie group $G$ can be identify with the space of left-invariant
vector fields equipped with the bracket operation. A
left-invariant sub-Riemannian structure on $G$ is simply a family
$\mathcal X=\{X_1,\dots,X_k\}$ of left-invariant vector fields
which generates the Lie algebra. The associated left-invariant
quadratic form on the cotangent bundle is given by (\ref{def-Q}).
We briefly illustrate this case by  examples of sub-Riemannian geometries
on  the Heisenberg group and on $\mbox{SU}(2)$.

\begin{example}
The Heisenberg group  $\mathbb H_3$ is $\mathbb R^3$ equipped
with the multiplication
$$(x,y,z)\cdot(x',y',z')=(x+x',y+y',z+z'+{\textstyle \frac{1}{2}}
(xy'-yx')).$$ It is easy to see that the left-invariant vector
fields equal to $d/dx, d/dy$ and $d/dz$ at $(0,0,0)$ are
$$X=d/dx- (y/2)d/dz, Y=d/dy+(x/2)d/dz \mbox{ and } Z=d/dz.$$
Moreover, $Z=[X,Y]$ and  $[X,Z]=[Y,Z]=0$. The ``canonical''
sub-Riemannian structure on $\mathbb H_3$ is associated with the
(minimal) H\"ormander family $\{X,Y\}$. This structure is
particularly adapted to $\mathbb H_3$ because it is homogeneous
with respect to the natural dilations $\delta_s(p)=(sx,sy,s^2z)$,
$p=(x,y,z)$,  that commutes with the group law. Namely, if
$d(p,p')$ is the sub-Riemannian distance associated with the
family $\{X,Y\}$, we have
$$d(\delta_s(p),\delta_s(p'))=s d(p,p'), \;\;p,p'\in \mathbb H_3.$$
No left-invariant Riemannian metrics can have this property. There
is an exact expression for the sub-Riemannian distance on $\mathbb
H_3$ (this is one of the very few cases of sub-Riemannian geometry
where such an exact formula exists). To connect with the notation
introduced in our discussion of Milnor frame on $3$-dimensional
Lie groups, observe that $f_1=Z/2$, $f_2=X$, $f_3=Y$ is a Milnor
frame for the left-invariant metric on $\mathbb H_3$ given at the
origin by $g_0=dx^2+dy^2+dz^2$. The sub-Riemannian structure on
$\mathbb H_3$ discussed above can be described as
$$Q = f_2\otimes f_2+ f_3\otimes f_3.$$
Note that this can be viewed as the limit of any family of
Riemannian metrics
$$g_t= A(t)f^1\otimes f^1+B(t)f^2\otimes f^2+C(t)f^3\otimes f^3,\;\;
t\in [0,T) $$
such that $\lim_T B=\lim_TC=1$ and $\lim_T A=\infty$. Indeed, in this case,
$$Q_t=A(t)^{-1}f_1\otimes f_1+B(t)^{-1}f_2\otimes f_2+C(t)^{-1}f_3\otimes f_3,
\;\; t\in [0,T)$$ which obviously tends to $Q$ (uniformly!). Note
that the existence of the dilations $\delta_s$, $s>0$, imply
immediately that the tangent cone at the identity element $e$ of
the pointed sub-Riemannian metric space $(\mathbb H_3, Q, e)$ is
that space itself.
\end{example}

\begin{example}
The group $\mbox{SU}(2)$ is the group of matrices
$$\left\{\left(\begin{array}{cc}a&b\\-\bar{b}&\bar{a}\end{array}\right):
a,b\in\mathbb C, |a|^2+|b|^2=1\right\}$$ which can also be
identified with the $3$-sphere $\mathbb S^3$. Its Lie algebra can
be identified with
$$\mathfrak{su}(2)= \left\{\left(\begin{array}{cc}i\alpha&\beta\\-\bar{\beta}&
i\alpha\end{array}\right):\alpha\in \mathbb R,\beta\in\mathbb
C\right\}.$$ Let $g_0$ be any left-invariant  Riemannian metric on
$\mbox{SU}(2)$ and $f_1,f_2,f_3$ be a Milnor frame as defined
earlier. Computation of the sectional curvatures (see below) show
that $(f_1,f_2,f_3)$ is always an orthonormal frame for the
standard round sphere metric on $\mbox{SU}(2)\simeq \mathbb S^3$.
Since $[f_i,f_j]=2f_k$ circularly, we can pick any two of these
vectors, say $f_2,f_3$, and consider the sub-Riemannian metrics
$$Q_{b,c}= bf_2\otimes f_2+cf_3\otimes f_3$$
where $b,c$ are fixed positive constants. These obviously appears as limits
of Riemannian metrics as in the case of the Heisenberg group discussed above.
The tangent cone of $\mbox{SU}(2)$ equipped with one of this sub-Riemannian
metric is, at any fixed point, the Heisenberg group equipped with
its canonical sub-Riemannian structure discussed above.
For a discussion of this
example and relation to the Hopf fibration, see \cite{Mont}.
\end{example}

\section{The cross curvature flow on the Heisenberg group}
Given a metric $g_0$ on the Heisenberg group (or on a $3$ manifold
of Heisenberg  type), we fix a Milnor frame $\{f_i\}_1^3$ such
that
$$[f_2,f_3]=2f_1, \;\;[f_3,f_1]=0,\;\;[f_1,f_2]=0$$
and (\ref{g_0})-(\ref{gflow}) hold.
Using \cite{milnor76}, the
sectional curvatures are:
$$
K(f_2 \wedge f_3) =  -\frac{3A}{BC}, \;\; K(f_3 \wedge f_1) =
K(f_1 \wedge f_2) =  \frac{A}{BC}.
$$
and the scalar curvature is $R = -2A/BC$.
The ODE for (+XCF) is given by
\begin{equation}
\left\{
\begin{aligned}
\frac{dA}{dt} =& 2\frac{A^3}{B^2C^2}, \\
\frac{dB}{dt} =&-6\frac{A^2}{BC^2},\\
\frac{dC}{dt} =& -6\frac{A^2}{B^2C}   .
\end{aligned}
\right .
\end{equation}
This case is very simple and admits a completely explicit
solution. Since the computations for (+XCF) are essentially the
same as for (-XCF), we refer the reader to \cite{cnsc1} for
details and simply write down the explicit maximal solution of
(-XCF) passing through $g_0$ at $t=0$.

\begin{theorem}\label{th-H1}
Given $A_0,B_0,C_0>0$ and $R_0=2A_0/(B_0C_0)$, set
$$T_0=T_b=1/(7R_0^2).$$
The maximal solution of {\em (-XCF)} through $g_0$ at $t=0$ is
defined on $(-T_0,\infty)$ and given by
$$\left\{\begin{array}{l}A(t)  =  A_0
(1+t/T_0)^{-\frac{1}{14}}\\
B(t) =  B_0 (1+t/T_0)^{\frac{3}{14}}\\
C(t) = C_0 (1+t/T_0)^{\frac{3}{14}}.
\end{array}\right.$$
The sectional curvatures are given by
$$-{\textstyle \frac{1}{3}}K(f_2 \wedge f_3)=K(f_1 \wedge f_2)
=K(f_3 \wedge f_1) = \frac{A_0}{B_0C_0}(1+t/T_0)^{-\frac{1}{2}}.$$
\end{theorem}
From the view point of the positive cross curvature flow (+XCF),
this theorem indicates that (+XCF) on the Heisenberg group
develops a singularity at the finite time $T_0=B_0^2C_0^2/(28
A_0^2)$. This singularity is of a type that is different from the
singularities usually observed in geometric flows which are dimensional
collapses such as pancake and cigar degeneracies.

Theorem \ref{th-H1} clearly indicates that it is natural to
re-scale the metric $g(t)$ by a factor of
$\psi(t)=(1+t/T_0)^{-3/14}$ (we will ignore the corresponding
change of time $\widetilde{t}=\int_0^t\psi(s)^2ds$ but note that
the backward blow-up time stays finite anyhow). Accordingly, we
set
$$\bar{g}(t)=
(1+t/T_0)^{-3/14} g(t).$$
\begin{theorem}\label{th-H2} Let $M$ be a complete locally homogeneous
$3$-manifolds of Heisenberg type with initial homogeneous metric
$g_0$ and associated Milnor frame $(f_1,f_2,f_3)$. Let $\bar{g}$
be as defined above.
\begin{enumerate}
\item If $M$ is compact, as $t$ tends to infinity,
the metric space $(M,g(t))$ converges in the
Gromov-Hausdorff sense to $\mathbb R^2$ with a flat metric (the
sectional curvature tends to $0$). If $M=\mathbb H_3$, as $t$ tends to
infinity,  the pointed metric space $(\mathbb H_3,e,g(t))$ converges in the
Gromov-Hausdorff sense to $\mathbb R^3$ with a flat metric.

\item As $t$ tends to $-T_0$, the
metric space $(M,\bar{g}(t))$ converges uniformly to the
sub-Riemannian metric space $(M, B_0^{-1}f_2\otimes
f_2+C^{-1}_0f_3\otimes f_3)$.
\end{enumerate}
\end{theorem}
\begin{remark}
In the first statement, the direction that collapses
is that of the center, i.e., $f_1$. It follows that to have Gromov-Hausdorff
convergence to a flat $\mathbb R^2$, it suffices to assume that $M$
is of the form $\mathbb H_3/\Gamma$ where $\Gamma$ is a discrete
subgroup of $\mathbb H_3$ with non-trivial intersection with the center.
If $\Gamma$ has trivial intersection with the center then
$M=\mathbb H_3/\Gamma$ is not compact and
converges in the pointed Gromov-Hausdorff
sense to a flat $\mathbb R^3$.
\end{remark}
\begin{remark}
In the second statement,
the identity map $\mbox{Id}$ is, in fact,
a bi-Lipschitz map between $(M,\bar{g}(t))$
and $(M, B_0^{-1}f_2\otimes
f_2+C^{-1}_0f_3\otimes f_3)$ and $\mbox{dil(Id)}$ defined at (\ref{def-dil})
tends to $0$ as $t$ tends
to $-T_0$.
\end{remark}
\begin{remark} \label{rem-heis}
Let $e_i^t$ be the unit vector for the metric $g(t)$
positively collinear to $f_i$. It is useful to look at the evolution of the
Lie algebra structure viewed from the perspective of the metric $g(t)$.
Namely, we have
$$[e^t_2,e^t_3]= \sqrt{A(t)/B(t)C(t)}e^t_1,
\;\;[e^t_1,e^t_2]=[e^t_3,e^t_1]=0.$$
As $t$ tends to $\infty$, $A(t)/B(t)C(t)$ tends to $0$ and the
non-trivial nilpotent structure converges to the trivial abelian structure on
$\mathbb R^3$. This is geometrically significant here since the exponential
map yields global coordinates.
\end{remark}

\section{The cross curvature flow on SU(2)}
Given a metric $g_0$ on $SU(2)$, we fix a Milnor frame such that
$[f_i,f_j]=2f_k$. For any metric $g=Af^1\otimes f^1+Bf^2\otimes
f^2+Cf^3\otimes f^3$, the sectional curvatures are (see, e.g.,
\cite [pg. 12]{chowknopf1}
\begin{align*}
K(f_2 \wedge f_3)&=\frac{(B-C)^2}{ABC}-\frac{3A}{BC}+\frac2B+\frac2C,\\
K(f_3 \wedge f_1)&=\frac{(C-A)^2}{ABC}-\frac{3B}{CA}+\frac2A+\frac2C,\\
K(f_1 \wedge
f_2)&=\frac{(A-B)^2}{ABC}-\frac{3C}{AB}+\frac2A+\frac2B.
\end{align*}

Together with the results obtained in \cite{cnsc1}, the asymptotic
behavior that will be proved in this section (see Theorem
\ref{th-su2-2} below) give the following full description of the
asymptotic behavior of the solution of the negative cross
curvature flow passing through a metric $g_0$ at time $0$.
\begin{theorem}\label{th-su2}
Let $g(t)=A(t)f^1\otimes f^1+B(t)f^2\otimes f^2+C(t)f^3\otimes
f^3$, $t\in (-T_b,T_f)$, be a maximal solution of the negative
cross curvature flow {\em (-XCF)} on $\mbox{\em SU}(2)$. Assume
that $A_0\ge B_0\ge C_0$ and set
$$\bar{g}(t) =\frac{B_0}{B(t)}g(t).$$
\begin{itemize} \item The time $T_f$ is finite. When
$t\rightarrow T_f$, $(\mbox{\em SU}(2), \bar{g}(t))$ converges
uniformly to the round sphere $\mathbb S^3_{\sqrt{B_0}}$.

\item If $A_0=B_0=C_0$, then $T_b=\infty$, $T_f=B_0^2/4$,
$B(t)=\sqrt{B_0^2-4t}$ and  $\bar{g}(t)=g_0$ (a round metric) for
all $t\in (-\infty,T_f)$ .

\item If $A_0=B_0>C_0$ (the Berger sphere case), then $T_b=\infty$
and, as  $t\rightarrow -\infty$, $(\mbox{\em SU}(2), \bar{g}(t))$
converges in the Gromov-Hausdorff sense  to a two-dimensional
round sphere (the sectional curvatures of $\bar{g}(t)$ containing
$f_3$ tends to $0$ whereas the one corresponding to $f_1\wedge
f_2$ tends to $4/B_0$).

\item If $A_0>B_0\ge C_0$ then $T_b$ is finite. As  $t\rightarrow
-T_b$, $(\mbox{\em SU}(2), \bar{g}(t))$ converges uniformly  to a
sub-Riemannian metric space $(\mbox{\em SU}(2), bf_2\otimes f_2+
cf_3\otimes f_3)$ with $b=B_0^{-1}$ and $c\in [B_0^{-1}, \infty)$
(if $B_0>C_0$ then $b>c$).
\end{itemize}
\end{theorem}
\begin{proof} Theorem 3 of \cite{cnsc1} gives $T_f\in (0,\infty)$
and the asymptotic $A,B,C\sim 2\sqrt{T-t}$ as $t$ tends to $T_f$.
The first statement follows. The other statements are consequences
of Theorem \ref{th-su2-2} below.
\end{proof}

\begin{remark}If one think about the global behavior of the
two-sided maximal flow lines of  the cross curvature flow in a
given Milnor frame on $\mbox{SU}(2)$ in terms of the value taken
at time $0$, $(A_0,B_0,C_0)$, in the first octant of $\mathbb
R^3$, one should distinguish $6$ regions, each  corresponding to a
strict order of the type $A>B>C$. These regions are preserved by
the flow and separated by planes of the type $A=B$ corresponding
to Berger sphere metrics. These planes are preserved by the flow.
The intersection of these planes is the line $A=B=C$ corresponding
to round metrics. This line is also preserved by the flow. In the
forward direction, the flow lines all approach the line $A=B=C$,
towards the point $(0,0,0)$. In the backward direction, in each of
the $6$ open regions, the largest component tends to infinity
whereas the two smaller components have finite distinct limits. On
the plane $A=B$, in the backward direction, $A=B$ tends to
infinity and $C$ has a finite limit.
\end{remark}

This section is devoted to the proof of the result that concern
the behavior at $-T_b$. For the rest of this section, we only
consider forward solutions of (+XCF).

From the sectional curvatures given above, we easily obtain the
ODEs corresponding to the positive cross curvature flow  (+XCF),
namely,
\begin{equation}\label{pdesu2}
\left \{
\begin{aligned}
\frac{dA}{dt}=&2\frac{AYZ}{(ABC)^2},\\
\frac{dB}{dt}=&2\frac{BZX}{(ABC)^2},\\
\frac{dC}{dt}=&2\frac{CXY}{(ABC)^2},
\end{aligned}
\right .
\end{equation}
where
\begin{align*}
X=&3A^2-(B-C)^2-2AB-2AC, \\
Y=&3B^2-(A-C)^2-2AB-2BC, \\
Z=&3C^2-(A-B)^2-2BC-2AC.
\end{align*}

Recall from \cite{cnsc1}, that we have
\begin{align}
\frac{d\ln(A/B)}{dt}=&\frac{8Z}{(ABC)^2}(B-A)(A+B-C),
\label{lnab}\\
\frac{d\ln(A/C)}{dt}=&\frac{8Y}{(ABC)^2}(A-C)(B-C-A), \label{lnac}\\
\frac{d\ln(B/C)}{dt}=&\frac{8X}{(ABC)^2}(B-C)(A-B-C), \label{lnbc}
\end{align}
and
\begin{align}
\frac{d(A-B)}{dt}=&\frac{2Z}{(ABC)^2}(B-A)[A^2+A(6B-2C)+(B-C)^2],\label{a-b}\\
\frac{d(B-C)}{dt}=&\frac{2X}{(ABC)^2}(C-B)[(A-B-C)^2+4BC],\label{b-c}\\
\frac{d(A-C)}{dt}=&\frac{2Y}{(ABC)^2}(C-A)((A-B)^2+6AC-2BC+C^2).
\label{a-c}
\end{align}

Without loss of generality we may assume that $A_0\ge B_0\ge C_0$
and it is easy to see from (\ref{a-b}) and (\ref{b-c}) that $A\ge
B\ge C$ is preserved along the flow. As a consequence, we have
\begin{align*}
Y&=(B-A)(A+B+2B-2C)-C^2\le -C^2<0,\\
Z&=-(A-B)^2+3C^2-2AC-2BC\le -C^2<0,
\end{align*}
and this implies the following Lemma.

\begin{lemma}\label{lem-mon}
Assume that $A_0\ge B_0\ge C_0$. Then $A$, $A/B$, $A/C$, $A-B$ and
$A-C$ are all nondecreasing along {\em (+XCF)}.
\end{lemma}

We will consider three cases.

{\bf Case 1: $A_0=B_0=C_0$.}

In this case, which is the round sphere,
$A(t)=B(t)=C(t)=\sqrt{A_0^2+4t}$, and the solution exists for all
time $t$.

{\bf Case 2: $A_0=B_0>C_0$.}

In this case, $A=B> C$ as long as the solution exist. The
equations simplify to
$$
\left \{
\begin{aligned}
\frac{dA}{dt}=&2\frac{4AC-3C^2}{A^3},\\
\frac{dC}{dt}=&2\frac{C^3}{A^4}.
\end{aligned}
\right .
$$

Clearly, $A$ and $C$ are increasing. Since $A\ddt A<8$, the
solution exists for all time $t\in [0,\infty)$.

\begin{lemma}Assume $A_0=B_0>C_0$. Then
$$\lim_{t \rightarrow \infty} A(t)=\infty~\text{and}~
\lim_{t \rightarrow \infty} C(t)=C_{\infty}<\infty.$$ Moreover, as
$t\rightarrow \infty$,
 $$A(t)\sim (24C_{\infty} t)^{1/3}~\text{and}~C(t)-
 C_{\infty}\sim - 2^{-3}3^{-1/3}C_{\infty}^{5/3}t^{-1/3}.$$
\end{lemma}

\begin{proof}
Assume that $\lim_{\infty} A=\eta<\infty.$ Since $A/C$ is
nondecreasing, we have
$$\ddt \ln(A/C) =\frac{8C(A-C)}{A^4}\ge \eta>0$$
and thus $\ln(A/C)>\eta t.$ This is a contradiction.

Now assume that $\lim_{\infty} C=\infty$. We first show that
$$\lim_{\infty} A/C=\infty.$$ Indeed, if
not, we must have $$\ddt \ln(A/C)\sim \eta/A^2,$$ for some
$\eta>0$. Hence $\int_0^{\infty} \frac{1}{A^2}<\infty.$ But we
also have
$$\ddt \ln C=\frac{2C^2}{A^4}\sim \eta \frac{1}{A^2}.$$
This contradicts the assumption that $\lim_{\infty} C=\infty$. So
$\lim_{\infty} A/C=\infty$ as desired.

We have $$\ddt \ln(A/C^4) =\frac{2C(4A-7C)}{A^4}.$$ Hence $A/C^4$
has a positive lower bound, say $\eta^2>0$. Now since
$$\ddt \frac{1}{A^{1/2}}\sim -4\frac{C}{A^{7/2}},$$ we have
$\int_0^{\infty} \frac{C}{A^{7/2}} <\infty.$ However $$\ddt
C=2\frac{C^2}{A^{1/2}}\times \frac{C}{A^{7/2}}$$ and
$\lim_{\infty} C=\infty$. It follows that
$$\int_0^{\infty} \frac{C}{A^{7/2}}> \frac{1}{\eta} \int_0^{\infty} \frac{C^2}{\sqrt{A}} \frac{C}{A^{7/2}}=  \infty. $$
This is a contradiction and we have proved that
$\lim_{\infty}C<\infty$.

 Now set $\lim_{\infty} C=C_{\infty}$. Then we
 easily check that
 $A(t)\sim (24C_{\infty} t)^{1/3},$ and $C(t)-C_\infty
 \sim -2^{-3}3^{-1/3}C_{\infty}^{5/3} t^{-1/3}$.
\end{proof}

{\bf Case 3: $A_0>B_0\ge C_0$.}

By Lemma \ref{lem-mon}, the condition $A_0>B_0$ implies that $A>B$
as long as the solution exits. Assume further that $A_0\geq 2B_0$.
Again, by Lemma \ref{lem-mon}, we have $A\geq 2B$ as long as the
solution exists.

\begin{lemma}\label{lem-mon2} Assume that $A_0\ge 2B_0$. Then $B/C$ is
nondecreasing. Furthermore $B$, $C$ and $B-C$ are non-increasing.
\end{lemma}
\begin{proof}
As $A\geq 2B$, we have $X>0$ and $A-B-C>0$. The desired
conclusions thus follow from (\ref{lnbc}), (\ref{b-c}) and
(\ref{pdesu2}).
\end{proof}

\begin{lemma} \label{lem-limit}Assume $A_0\ge 2B_0$.
There exists a $T<\infty$, such that $\lim_{t\rightarrow T}
C(t)=0$, $\lim_{t\rightarrow T} B(t)=0$ and $\lim_{t\rightarrow T}
A(t)=\infty$.
\end{lemma}

\begin{proof}
As $A\geq 2B$, we have \begin{equation}\label{X<} \frac34
A^2<X<3A^2.\end{equation}
 It follows that
$$-\ddt B=2\frac{B|Z|X}{A^2B^2C^2}> \frac{3}{2B},$$ or $$-\ddt
B^2>3.$$ Hence the solution can only exist up to some finite time
$T$, i.e., there exists $T<\infty$, such that either $\lim_{T}
C=0$ or $\lim_{T} A=\infty$.

\begin{claim}\label{3.1} If $\lim_{t\rightarrow T} C(t)=0$, then
$\lim_{t\rightarrow T} B(t)=0$ and $\lim_{t\rightarrow T}
A(t)=\infty$.
\end{claim}

To prove the claim, suppose $\lim_{T} C=0$, $\lim_{T} B=B(T)>0$
and $\lim_{T} A=A(T)<\infty$. Since
$$\ddt C=2\frac{CXY}{A^2B^2C^2},$$ and $$0<\lim_T
\frac{X|Y|}{A^2B^2}=\lim_T \frac{X}{A^2}
\frac{(A-B)(A+3B-2C)}{B^2}=\eta_1<\infty,$$ we have
$$-\ddt C^2 \sim 4\eta_1,$$ and thus $C^2(t)\sim 4\eta (T-t).$
As $t\rightarrow T$,  $|Z|\sim (A-B)^2$ has a positive finite
limit. By (\ref{X<}), $X$ also has a positive finite limit. Hence
$$-\ddt B^2=\frac{4|Z|X}{A^2C^2}\sim \frac{\eta_2}{T-t}.$$
This contradicts $\lim _TB=B(T)>0$.

Now, we either have $\lim_{T} B=0$ or $\lim_{T} A=\infty$ and, in
particular, $\lim_{T} \frac{A}{B}=\infty$. Further, as
$t\rightarrow T$, we have
$$X\sim 3A^2,\;\;Y\sim -A^2,\; \mbox{ and }Z\sim -A^2.$$

Assume that $\lim_{T} B(t)=B(T)>0$. Observe that $$-\ddt \ln
(B-C)=\frac{2X}{A^2B^2C^2}[(A-B-C)^2+4BC]\sim \frac{6
A^2}{B^2C^2},$$ hence $$\int_0^T \frac{A^2}{B^2C^2}<\infty.$$ But
we also have $$\ddt \ln A=\frac{2YZ}{A^2B^2C^2}\sim
\frac{2A^2}{B^2C^2}.$$ As $\lim_{T} A=\infty$, this is a
contradiction. It follows that we must have $\lim_{T} B=0$. A
similar argument show that $\lim_{T} A=\infty$. This finishes the
proof of the claim \ref{3.1}.

To finish the proof of Lemma \ref{lem-limit}, it suffices to rule
out the case $\lim_{T} C>0$. Assume $\lim_{T} C=C(T)>0$, then
$\lim_{T} B=B(T)>0$ and $\lim_{T} A=\infty$. So, as $t\rightarrow
T$, we have $X\sim 3A^2$, $Y\sim -A^2$ and $Z\sim -A^2$. Hence
$$\ddt A=\frac{2AYZ}{A^2B^2C^2} \sim \eta A^3$$ with
$\eta=2/(B(T)C(T))^2$. It follows that
$$A(t)\sim \frac{1}{\sqrt{2\eta(T-t)}}.$$ We also have
$$\ddt B^2\sim -12 \frac{A^2}{C(T)^2} \sim -\frac{3B(T)^2}{T-t}.$$
This contradicts $B(T)>0$ and thus, we must have $\lim_TC=0$ as
desired.
\end{proof}

\begin{lemma}
Assume $A_0\geq 2B_0$. As $t\rightarrow T$, there exist positive
finite constants $\eta_1$, $\eta_2$ and $\eta_3$ such that
\begin{equation}\label{3eq}
\left \{
\begin{aligned}
A \sim & \eta_1 (T-t)^{-\frac{1}{14}},\\
B \sim & \eta_2 (T-t)^{\frac{3}{14}},\\
C \sim & \eta_3 (T-t)^{\frac{3}{14}}.
\end{aligned}
\right .
\end{equation}
\end{lemma}
\begin{proof} As a first step in the proof of this lemma we
show the following.
\begin{claim}\label{3.2}
As $t\rightarrow T$, we have $\lim_{t\rightarrow T}
\frac{B}{C}=\eta \in[1,\infty)$.
\end{claim}
 We have
$$
\left \{
\begin{aligned}
\frac{dA}{dt}\sim &2\frac{A^3}{(BC)^2},\\
\frac{dB}{dt}\sim &-6\frac{A^2}{BC^2},\\
\frac{dC}{dt}\sim &-6\frac{A^2}{B^2C}.
\end{aligned}
\right .
$$
Moreover, by (\ref{lnbc}), there exists $\eta_1\in (0,\infty)$
such that
$$\ddt \ln (B/C)=\frac{8X}{A^2BC^2} \frac{B-C}{B} (A-B-C)
\sim \eta_1 \frac{A}{BC^2}.$$ Here we used that $1-\frac{C}{B}$ is
increasing under the flow. We also have $-\ddt B\sim
6\frac{A^2}{BC^2},$ so $\int_0^T \frac{A^2}{BC^2} <\infty.$  As
$A\ge A_0$ this implies that $\int_0^T \frac{A}{BC^2}<\infty.$
Hence $\lim_{ T} \frac{B}{C}=\eta <\infty$ as claimed.

\begin{claim}
\label{a3b} As $t\rightarrow T$, we have
 $\lim_{t\rightarrow T} A^3B=\eta_1$ and
$\lim_{t\rightarrow T} A^3C=\eta_2$ with $\eta_1,\eta_2\in
(0,\infty)$.
\end{claim}

To prove this claim, we compute
$$\ddt (A^{\alpha}B)=\frac{2A^{\alpha}BZ}{A^2B^2C^2}(X+\alpha
Y).
$$
In particular, if $\alpha=3$, $$\ddt
(A^{3}B)=\frac{8A^{3}BZ}{A^2B^2C^2}[2B^2-C^2-BC+AC-2AB]>0.$$ In
general,
$$X+\alpha Y
=(3-\alpha)A^2-(B-C)^2-2AB-2AC+3\alpha B^2-\alpha C^2 +2\alpha
AC-2\alpha AB-2\alpha BC.$$ In particular, if $\alpha\neq3$, then
as $t\rightarrow T$,
$$
X+\alpha Y \sim (3-\alpha)A^2.$$

For $\alpha=2$, we see that $\ddt (A^{2}B)<0$ for $t$ close to $T$
because $Z<0$. Hence $A^2B$ is bounded from above. Now, Claim
\ref{3.2} and the above formula yields
$$\ddt
(A^{3}B)\sim 8(2-1/\eta) \frac{A^4}{C^2}.$$ Since $\ddt B\sim -6
\frac{A^2}{BC^2},$ and $\lim B=0$, we have $\int_0^T
\frac{A^2}{BC^2} <\infty.$ As $A^2B$ is bounded from above and
$A^4/C^2= (A^2B)\times A^2/(BC^2)$, it follows that $\int_0^T
\frac{A^4}{C^2}<\infty.$ Hence
$$\lim_{ T} A^3B=\eta_1\in (0,\infty).$$ As $\lim
B/C=\eta$, this also yields $\lim_T A^3C= \eta_1/\eta$ as claimed.

Now using Claim \ref{a3b} and (\ref{3eq}), we obtain
$$
\frac{dA}{dt}\sim  2\eta_1^{-2}\eta_2^{-2} A^{15}.
$$
Hence,
$$
\left \{
\begin{aligned}
A \sim & \eta_3 (T-t)^{-\frac{1}{14}},\\
B \sim & \eta_4 (T-t)^{\frac{3}{14}},\\
C \sim & \eta_5 (T-t)^{\frac{3}{14}}.
\end{aligned}
\right .
$$
\end{proof}

In order to finish the study of the behavior of (+XCF) on
$\mbox{SU}(2)$, we are left to show that if $A_0>B_0\ge C_0$ then
there exists a time $t_0\ge 0$ such that $A(t_0)\ge 2B(t_0)$. We
start with the following lemma.

\begin{lemma}
If there exist a time $t_0$ such that $X(t_0)>0$, then for all
time $t\ge t_0$, we have
$$\frac{X}{A^2}\ge \frac{X(t_0)}{A(t_0)}=\eta>0.$$ Moreover, there
exists a time $t_0'$ such that $A(t'_0)\ge 2B(t_0')$.
\end{lemma}

\begin{proof}
We first show that $\frac{X}{A^2}>0$ for $t>t_0$. Let $t_1>t_0$ be
the first time such that $X(t_1)=0$ (if it exists).
 Recall that $\frac{B}{A}$ and $\frac{C}{A}$ are decreasing.
 Moreover, for all
$t\in [t_0,t_1)$, we have $$\ddt \ln
(B-C)=-\frac{2X}{A^2B^2C^2}[(A-B-C)^2+4BC]<0,$$ hence
$\frac{B-C}{A}$ is decreasing as well on this interval.
 Since
$$\frac{X}{A^2}=3-\frac{(B-C)^2}{A^2}-\frac{2B}{A}-\frac{2C}{A},$$
$X/A^2$ is increasing on $[t_0, t_1)$. This contradicts
$X(t_1)=0$, hence we have $\frac{X}{A^2}>0$, for all $t>t_0$. When
$X$ is positive, $B-C$ is decreasing and thus  $X/A^2$ is
increasing. It follows that  $\frac{X}{A^2}>\eta>0$.

Now, assume that for all time $t>t_0$, $A(t)<2B(t)$. Since $A$ is
increasing and $B$ and $C$ are decreasing, we have $A(t_0)\leq
A(t) <2B(t) \leq 2B(t_0)<\infty$, for $t>t_0$. Since
$$-\ddt B=\frac{2B|Z|X}{A^2B^2C^2}\geq \frac{\eta'}{B},$$
for some $\eta'>0$, the solution can only exist up to some finite
time $T<\infty$. This means that $\lim_{T} C=0$. Since
$\lim_TA=A(T)$ and $\lim_TB=B(T)$ are positive and finite,  as
$t\rightarrow T$, there is $\eta_1>0$ such that
$$\ddt C^2 \sim -\eta_1.$$
Hence $C\sim \sqrt{\eta_1(T-t)}.$ This implies $$\ddt B^2 \sim
-\frac{\eta_2}{T-t}.$$ This contradicts the fact that $B^2\ge 0$.

\end{proof}

Our final task is to show that it is not possible that
 $A<2B$ and $X\leq 0$ for all $t$. To this end, assume that
 $A<2B$ and $X\le 0$ on the interval $[0,T)$ on which the solution
 exists ($T$ might be $\infty$).
 Then $A$, $B$, $A/B$ and $C$ are all non-decreasing on $[0,T)$.
We claim that
$$\sup_{[0,T)} \left\{B/C\right\}=\eta <\infty.$$
Otherwise, since $B<A$ and $B/A$ is non-increasing, we have
$$\limsup_{t\rightarrow T} \frac{X}{A^2}= \lim_{t\rightarrow T}
(3-\frac{2B}{A}-\frac{B^2}{A^2}+\frac{-C^2+2BC-2AC}{A^2})>0,$$
which contradicts $X\leq 0$.

Hence we have $A<2B<2\eta C$. Observe that $Y,Z$ are negative and
$$|Y|=A^2+C^2-2AC+2AB+2BC-3B^2<A^2+C^2+2AB+2BC<11B^2,$$ and
$$|Z|=A^2+B^2-2AB+2BC+2AC-3C^2<A^2+B^2+2BC+2AC<\eta_1 C^2,$$
for some $\eta_1\in (0,\infty)$. It follows that
$$\ddt A=\frac{2AYZ}{A^2B^2C^2}<\eta_2 \frac{1}{A},$$
for some $\eta_2\in (0,\infty)$. This implies that $T=\infty$,
i.e., the solution exists for all time. We also have
$$|X|=B^2+C^2+2AB+2AC-2BC-3A^2<6A^2. $$ Hence
$$\ddt B=\frac{2BZX}{A^2B^2C^2}<\eta_3 \frac{1}{B},$$
for some $\eta_3\in (0,\infty)$ and $B^2\leq B_0^2+2\eta_3 t.$
Now, since
$$\ddt \ln (A/B)=\frac{8|Z|}{A^2B^2C^2}
(A-B)(A+B-C)\ge \frac{8}{B^2}(1-\frac{B}{A})\ge
\frac{\eta_4}{B^2},$$ we have
$$\ln (A/B) -\ln
(A_0/B_0)\ge \eta_4 \int_0^{t} \frac{1}{B^2}\ge\eta_4 \int_0^{t}
\frac{ds}{B_0^2+2\eta_3 s}.$$ The right hand side goes to $\infty$
as $t\rightarrow \infty$, while we have
$$\lim_{t\rightarrow \infty} \ln (A/B)\le \ln 2.$$ This is
the desired contradiction. The conclusion is that, if $A_0>B_0\ge
C_0$ then there must exists a time $t_0$, such that $A(t_0)\geq
2B(t_0)$.

\begin{theorem}\label{th-su2-2}
On $\mbox{\em SU}(2)$, for the positive cross curvature flow {\em
(+XCF)} and any choice of initial data $A_0 \geq B_0 \geq C_0>0$,
we have:
\begin{itemize}
\item If $A_0=B_0=C_0$, then the solution  exists for all time
$t$, and we have $$A(t)=B(t)=C(t)=\sqrt{A_0^2+4t}.$$

\item If $A_0=B_0>C_0$, then the solution exists for all time $t$
and $\lim_{t \rightarrow \infty} C(t)=C_\infty\in (0,\infty)$.
Moreover, as $t$ tends to $\infty$, we have
 $$A(t)=B(t) \sim (24C_\infty t)^{1/3},~C(t)-C_\infty\sim
 -2^{-3}3^{-1/3}C_\infty^{5/3} t^{-1/3}.$$

\item if $A_0>B_0\ge C_0$, then there exists a finite time $T>0$,
such that the solution exists on $[0,T)$ and, as $t\to T$,
$$
A \sim  \eta_1 (T-t)^{-\frac{1}{14}},~ B \sim  \eta_2
(T-t)^{\frac{3}{14}},~and~ C \sim  \eta_3 (T-t)^{\frac{3}{14}},
$$
for some finite positive constants $\eta_1,\eta_2,\eta_3$.
\end{itemize}
\end{theorem}

\section{The cross curvature flow on $E(1,1)$ (Sol geometry)}
A model for the geometry $E(1,1)$ is the group $\mathbb R\ltimes
\mathbb R^2$ where the action of $\mathbb R$ on $\mathbb R^2$ is
given by
$\left(\begin{array}{cc}e^x&0\\0&e^{-x}\end{array}\right)$. This
group is sometimes called Sol. In other words, Sol is $\mathbb
R^3$ with the multiplication law
$$(x,y,z)\cdot (x',y',z')= (x+x',y+e^xy',z+e^{-x}z').$$
If we denote by $X,Y,Z$ the left-invariant vector fields equal to
$\partial/\partial x$, $\partial /\partial y$ and $\partial
/\partial z$ at $(0,0,0)$ then $X=\partial/\partial x$,
$Y=e^x\partial/\partial y$ and $Z=e^{-x}\partial/\partial z$.
Hence, $[Y,Z]=0$, $[X,Y]=Y$, $[X,Z]=-Z$. For any left-invariant
metric equal to $a(dx)^2+b(dy)^2+c(dz)^2$ at $(0,0,0)$, a Milnor
frame is $f_1= Y+Z$, $f_2=2X$, $f_3=Y-Z$. Conversely, given a metric $g_0$
and a Milnor frame for that metric, we can define elements $X,Y,Z$
of the Lie algebra by $X=(1/2)f_2$, $Y=(1/2)(f_1+f_3)$ and $Z=(1/2)(f_1-f_3)$.
In the coordinate system induced by the exponential map and the basis
$X,Y,Z$, the group law has the form given above.
It is useful to observe that, for any $r>0$, the metrics
$$g=af^1\otimes f^1+bf^2\otimes f^2+cf^3\otimes f^3
\mbox{ and } g= raf^1\otimes f^1+bf^2\otimes f^2+rcf^3\otimes
f^3$$ yields isometric manifolds, the isometry being induced by
the Lie algebra isomorphism  $f_1\mapsto \sqrt{r}f_1$, $f_2\mapsto
f_2$, $f_3\mapsto \sqrt{r} f_3$ and its inverse.

For later purpose, we introduce the manifold
$$E_0(1,1)=E(1,1)/\mathbb X$$
where $\mathbb X$ stands for the copy of the integers
$\{(k,0,0): k\in \mathbb Z\}$ sitting in $E(1,1)$. Any left-invariant
metric on $E(1,1)$ induces a locally homogeneous Riemannian structure on
$E_0(1,1)$. Moreover, the remark made above concerning isometric structures
on $E(1,1)$ is also valid for $E_0(1,1)$.
We will be particularly interested in those metric $g_\beta$, $\beta>0$,
of the form
\begin{equation}\label{def-gc}
g_\beta=\beta f^2\otimes f^2 +
f^1\otimes f^1 +f^3\otimes f^3.
\end{equation}

Given a left-invariant metric $g_0$, we fix a
Milnor frame $\{f_i\}_1^3$ such that
$$[f_2,f_3]=2f_1, \;\;[f_3,f_1]=0,\;\;[f_1,f_2]=-2f_3.$$

The sectional curvatures are:
\[
K(f_2 \wedge f_3) = \frac{(A-C)^2-4A^2}{ABC},
\]

\[
K(f_3 \wedge f_1) = \frac{(A+C)^2}{ABC},
\]

\[
K(f_1 \wedge f_2) = \frac{(A-C)^2-4C^2}{ABC}.
\]
We first state the theorem describing the forward and backward
asymptotic behaviors of (-XCF) on $E(1,1)$. Because of the
symmetry between $f_1$ and $f_3$, one can assume that $A\ge C$.

\begin{theorem}\label{th-E11}
Let $g(t)=A(t)f^1\otimes f^1+B(t)f^2\otimes f^2+C(t)f^3\otimes
f^3$, $t\in (-T_b,T_f)$, be a maximal solution of the negative
cross curvature flow {\em (-XCF)} on a complete locally
homogeneous manifold of type $E(1,1)$. Assume that $A_0\ge C_0$
and set
$$\bar{g}(t) =\frac{B_0}{B(t)}g(t).$$
\begin{itemize} \item The time $T_f$ is finite. If $M$ is compact, as
$t\rightarrow T_f$, the metric space $(M, \bar{g}(t))$ converges
in the Gromov-Hausdorff sense to $(E_0(1,1),g_\beta)$ with
$g_\beta$ as in (\ref{def-gc}) for some $\beta>0$. If $M=E(1,1)$,
then  the metric space $(E(1,1), \bar{g}(t))$ converges in the
Gromov-Hausdorff sense to $(E(1,1),g_\beta)$ for some $\beta>0$.
\item If $A_0=C_0$, then $T_b=\infty$ and $A(t)=C(t)=A_0B_0/B(t),$
$B(t)=\sqrt{B_0^2-64t}$. If $M$ is compact, as $t\rightarrow
-T_b=-\infty$, the metric space $(M,\bar{g}(t))$ (resp.
$(M,g(t))$)   converges in the Gromov-Hausdorff sense to a circle
(resp. to a line).

\item If $A_0>C_0$, then $T_b<\infty$ and, as  $t\rightarrow
-T_b$, the metric space $(M, \bar{g}(t))$ converges uniformly to
the sub-Riemannian metric space $(M, bf_2\otimes f_2+ cf_3\otimes
f_3)$ with $b,c\in (0, \infty)$.
\end{itemize}
\end{theorem}
\begin{proof} Theorem 2 from \cite{cnsc1} gives that $T_f$ is finite and
$A,C\sim E_1/\sqrt{T_f-t}$, $A-C\sim E_2\sqrt{T_f-t}$, $B\sim
\sqrt{64(T_f-t)}$, as $t$ tends to $T_f$. The desired convergence
follows. For the behavior when $t$ tends to $-T_b$, see Theorem
\ref{th-E11-2} below.
\end{proof}

The rest of this section is devoted to the asymptotic behavior of
the flow ($*$XCF) on $E(1,1)$. From the sectional curvature given
above, we obtain the equations
\begin{equation}
\left \{
\begin{aligned}
\frac{dA}{dt} =&\frac{2AYZ}{(ABC)^2}
\\
\frac{dB}{dt} =&\frac{2BZX}{(ABC)^2} 
\\
\frac{dC}{dt} =&\frac{2CXY}{(ABC)^2}
\end{aligned}
\right .
\end{equation}
where
\begin{align*}
X=&(A+C)(3A-C), \\
Y=&-(A+C)^2, \\
Z=&-(A+C)(A-3C).
\end{align*}

First, consider the case when  $A_0=C_0$. Then $A(t) = C(t) $ as
long as the solution exists and
\[
\frac{dB}{dt} = \frac{32}{B}.
\]
Hence $B=\sqrt{B_0^2+64t}$. It follows that
\[
\frac{d \ln A}{dt} = -\frac{32}{B_0^2+64t},
\]
which gives
\[
A(t)=C(t) = \frac{A_0B_0}{ \sqrt{B_0^2+64t}}.
\]

Second, we assume that  $A_0\ne C_0$. Because of the symmetry
between $f_1$ and $f_3$, we may assume without loss of generality
that $A_0 > C_0$. It  immediately follows that $C$ is decreasing.
\begin{lemma}Assume that $A_0>C_0$. Then $C$ is decreasing,
$A-C$, $A/C$ and $A-3C$ are increasing as long as the solution
exists. In particular, we have $A>C$.
\end{lemma}
\begin{proof} Observe that
\begin{align*}
&\frac{d(A-C)}{dt}=2\frac{(A+C)^4}{(ABC)^2}(A-C),\\
&\frac{d\ln(A/C)}{dt}=8\frac{(A+C)^3}{(ABC)^2}(A-C),\\
&\frac{d(A-3C)}{dt}=2\frac{(A+C)^3}{(ABC)^2}(A^2+6AC-3C^2).
\end{align*}
The stated results follow.\end{proof}

\begin{lemma}\label{lem-a3c}
Assume that $A_0>3C_0$. Then $A$ is increasing,
$B$ is decreasing. Moreover, there exists a finite time $T$ such
that
$$\lim_{t\rightarrow T}A(t)=\infty,\;\;
\lim_{t\rightarrow T}B(t)=0,\;\; \lim_{t\rightarrow T}C(t)=0.$$
\end{lemma}
\begin{proof}
The previous lemma shows that $A/C>A_0/C_0>3$ as long as the
solution exists. The monotonicity of $A$ and $B$ immediately
follows. Further, we have
$$-B \ddt B=2\frac{(3A-C)(A-3C)(A+C)^2}{A^2C^2}
>4(1-\frac{3C}{A})(1+\frac{A}{C})^2>\eta.$$
This implies that the solution can only exist up to some finite
time $T$. As
$$\ddt \ln (B/C)=\frac{8(A+C)^2(3A-C)}{A^2B^2C}>0,$$
$B/C$ is increasing. Hence if $\lim_{t\rightarrow T} B(t)=0$, then
$\lim_{t\rightarrow T} C(t)=0$ as well.

Assume that we have $\lim_{T} A=A(T)<\infty$ and $\lim_{T}
B=B(T)>0$. Then we must have $\lim_{T} C=0$ (otherwise, the
solution could be extended past $T$). As
$$\ddt C^2=\frac{4XY}{(AB)^2},$$ it then follows that
$C(t)^2\sim \eta_1 (T-t)$, for some positive finite $\eta_1$. This
implies
$$\ddt A=\frac{2AYZ}{(ABC)^2} \sim \frac{\eta_2}{T-t}$$
and thus $A\rightarrow \infty$. This is a contradiction.

Assume next that $\lim_{T} A=A(T)<\infty$ and $\lim_{T} B=0$. Then
$\lim_TC=0$ and
 $$B^2 \ddt
C^2=-\frac{4(A+C)^3(3A-C)}{A^2} \sim -12A^2,$$ $$C^2 \ddt
B^2=\frac{4(3A-C)(3C-A)(A+C)^2}{A^2} \sim -12A^2.$$ Hence $$\ddt
(B^2C^2)\sim -24A^2$$ and $B^2C^2 \sim \eta_3(T-t)$ with
$\eta_3\in (0,\infty)$. Plugging this into the formula for $\ddt
A$ shows that $A \rightarrow \infty$, which is a contradiction.

This shows that, as stated in Lemma \ref{lem-a3c}
$\lim_TA=\infty$. To see that $\lim_T B=0$, we compute
\begin{equation}\label{eq.alpha}
\ddt (A^{\alpha}B)=\frac{2A^{\alpha}BZ}{A^2B^2C^2}(X+\alpha Y).
\end{equation} If $0<\alpha <3$, then as $t\rightarrow T$,
$$X+\alpha Y =(A+C)[(3-\alpha)A-(1+\alpha )C]>0.$$ As $Z$ is negative,
this implies that $\ddt A^\alpha B<0$. Thus $A^{\alpha}B$ is
decreasing. Since $A\rightarrow \infty$, we also have
$B\rightarrow 0$ and thus $C\rightarrow 0$.
\end{proof}

\begin{lemma}\label{lem-a3c'}
Assume that $A_0>3C_0$ and let $T$ be as in Lemma
\ref{lem-a3c}. Then there are positive finite constants $\eta_i$,
$i=1,2,3$ such that, as $t$ tends to $T$,
$$
A \sim \eta_1 (T-t)^{-\frac{1}{14}},\;\;
B \sim  \eta_2 (T-t)^{\frac{3}{14}},\;\;
C \sim  \eta_3 (T-t)^{\frac{3}{14}}.
$$
\end{lemma}
\begin{proof}
Taking $\alpha=3$ in (\ref{eq.alpha}) yields $$\ddt
(A^{3}B)=\frac{2AZ}{BC^2}(3Y+X)=-\frac{8AZ}{BC}(A+C)>0.$$
Similarly,
$$\ddt
(A^{3}C)=\frac{2AY}{B^2C}(3Z+X)=-\frac{16A}{B^2}(A+C)^3<0.$$
 Now,  the  technique used in the proof of Claim \ref{a3b}
shows that
$$\lim_{t\rightarrow T} A^3B<\infty\mbox{ and }
\lim_{t\rightarrow T} A^3C>0.$$ As $t$ tends to $T$, this yields
$$
\frac{dA}{dt}\sim  \eta_1 A^{15},\;\;
\frac{dB}{dt}\sim  -\eta_2 B^{-\frac{11}{3}},\;\;
\frac{dC}{dt}\sim -\eta_3 C^{-\frac{11}{3}},
$$
for positive finite $\eta_i$, $i=1,2,3$. The asymptotic stated in
Lemma \ref{lem-a3c'} follows. \end{proof}

We are left with the task to rule out the possibility that
$C<A\leq 3C$ as long as the solution exists. Assume that $C<A\leq
3C$. Then $A$ and $C$ are decreasing and $B$ is increasing. Since
we have $$B\ddt B=2\frac{(3A-C)(3C-A)(A+C)^2}{(AC)^2}<\eta,$$ we
have
$$B\leq \sqrt{B_0^2+2\eta t}.$$
It follows that
$$\ddt \ln (A/C) =8 \frac{(A+C)^3}{A^2C} \frac{1}{B^2}
(\frac{A}{C} -1) > \frac{\eta}{B_0^2+2\eta t}.$$ As $\ln
(A/C) \leq \ln 3$,  the solution can only exist up to some
finite time $T$, i.e., there exist $T<\infty$, such that
$\lim_{t\rightarrow T} A(t)=\lim_{t\rightarrow T} C(t)=0$. Hence
$\lim_{t\rightarrow T} (A-C)(t)=0$. This contradicts the fact that
$A-C$ is increasing.  It follows that there must exist a finite
time $t_0$, such that $A(t_0)>3C(t_0)$. We have proved the
following theorem.

\begin{theorem}\label{th-E11-2}
On $E(1,1)$, for any given initial data $A_0,B_0,C_0>0$ with
$A_0\ge C_0$, the positive cross curvature flow behaves as
follows.
\begin{itemize}
\item If $A_0=C_0$, then the solution of {\em (+XCF)} exists on
$[0,\infty)$ and is given by
$$A(t)=C(t)=\frac{A_0B_0}{\sqrt{B_0^2+64t}},\;\;B(t)=\sqrt{B_0^2+64t}.$$

\item If $A_0>C_0$, then there exists a positive finite time $T_b$
such that the solution of {\em (+XCF)} exists on $[0,T_b)$.
Moreover, as $t\to T_b$,
$$
A \sim  \eta_1 (T-t)^{-\frac{1}{14}},~ B \sim  \eta_2
(T-t)^{\frac{3}{14}},~and~ C \sim  \eta_3 (T-t)^{\frac{3}{14}}
$$
with $\eta_1,\eta_2,\eta_3\in (0,\infty)$.
\end{itemize}
\end{theorem}

\section{The cross curvature flow on $E(2)$}
Recall that a realization of $\widetilde{E(2)}$ is $\mathbb R^2
\rtimes \mathbb R$ where the action of $\mathbb R$ on $\mathbb
R^2$ is by rotation. Namely, $\widetilde{E(2)}$ is $\mathbb R^3$
equipped with the product
$$(x,y,z)\cdot (x',y',z')=( x+ x'\cos z +y'\sin z, y+ x'\cos
z-y'\sin z,z+z').$$ Note that the left-invariant metrics on
$\widetilde{E(2)}$ equal to $a[(dx)^2+(dy)^2]+c(dz)^2$ at the
identity element are actually flat metrics on $\mathbb R^3$.

Given a left-invariant metric $g_0$, we fix a Milnor frame
$\{f_i\}_1^3$ such that
$$[f_2,f_3]=2f_1, \;\;[f_3,f_1]=2f_2,\;\;[f_1,f_2]=0.$$
These relations imply that, in the global coordinate introduce
above, $f_1,f_2$ are in the span of $\partial/\partial x,
\partial /\partial y$ whereas $f_3$ contains a $\partial /\partial z$
component.

The sectional curvatures are:
\[
K(f_2 \wedge f_3) = \frac{1}{ABC} (B-A)(B+3A),
\]

\[
K(f_3 \wedge f_1) = \frac{1}{ABC} (A-B)(A+3B),
\]

\[
K(f_1 \wedge f_2) = \frac{1}{ABC} (A-B)^2.
\]
As in previous sections, we first state a theorem that describes
both the forward and backward behavior of the negative cross
curvature flow on $E(2)$. The forward behavior was studied in
\cite[sect. 6]{cnsc1}. The backward behavior is studied below in
the form of the forward behavior of the positive cross curvature
flow. Because of the symmetry between $f_1$ and $f_2$ we can assume
that $A_0\ge B_0$.
\begin{theorem}Let $g(t)=A(t)f^1\otimes f^1+B(t)f^2\otimes  f^2  +
C(t)f^3\otimes f^3$, $t\in (-T_b,T_f)$, be a maximal solution of
the negative cross curvature flow on a complete locally
homogeneous manifold of type  $E(2)$. Set
$$\bar{g}(t)=\frac{C_0}{C(t)} g(t).$$
\begin{itemize}
\item  Assume that $A_0=B_0$. Then $T_b=T_f=\infty$,
$A(t)=B(t)=A_0$ and $C(t)=C_0$ for all $t$. In fact, as a
Riemannian manifold, $(\widetilde{E(2)},g)$ is $\mathbb R^3$
equipped with a flat metric. \item Assume that $A_0>B_0$.  Then
$T_f=\infty$ and $T_b<\infty$.
\begin{itemize}
\item If $M=\widetilde{E(2)}$, as $t$ tends to $T_f=\infty$,
$(\widetilde{E(2)}, g(t))$ converges in the Gromov-Hausdorff sense
to $\mathbb R^3$ (the sectional curvatures go to $zero$). If $M$
is compact, $(M,g(t))$ does not converge  to a metric space as $t$
tends to infinity but $(M,\bar{g}(t))$ converge in the
Gromov-Hausdorff sense to a circle. \item In all cases, as $t$
tends to $-T_b$, $(M,\bar{g}(t))$ converges uniformly to the
sub-Riemannian metric space $(E(2), bf_2\otimes f_2+ cf_3\otimes
f_3)$ for some $b,c\in (0,\infty)$.
\end{itemize}
\end{itemize}
\end{theorem}
\begin{proof}The first case is trivial. When $A_0>B_0$, the results in
the forward direction follow from Theorem 5 of \cite{cnsc1} which
gives $T_f=\infty$ and $A\sim E_1+E_2t^{-1/6}$, $B\sim
E_1-E_2t^{-1/6}$ and $C\sim (8\sqrt{6}E_2/E_1)t^{1/3}$, as $t$
tends to infinity.

The result in the backward direction follow from Theorem \ref{th-E22} below.
\end{proof}
\begin{remark} For $t>0$, in the
orthonormal frame $(e^t_i)1^3$ on $(\widetilde{E(2)},g(t))$, with
$e^t_i$ positively collinear to $f_i$, we have $[e^t_2,e^t_3]\sim
\delta t^{-1/6} e^t_1$, $[e^t_3,e^t_1]\sim \delta t^{-1/6} e^t_2$,
$[e_1^t,e_2^t]=0$ for some $\delta>0$. In other words, the group
structure viewed in this frame tends to the abelian structure of
$\mathbb R^3$. This is similar to what happens on the Heisenberg
group. Compare Remark \ref{rem-heis}.
\end{remark}

The rest of this section is devoted to the study of the positive
cross curvature (+XCF) on $E(2)$. Hence $g(t)=A(t)f^1\otimes
f^1+B(t)f^2\otimes f^2+C(t)f^3\otimes f^3$ is a solution of (+XCF)
and thus
 $A,B,C$ satisfy the following
equations
\begin{equation}\label{pdee2}
\left\{
\begin{aligned}
\frac{dA}{dt}=&\frac{2AYZ}{(ABC)^2},\\
\frac{dB}{dt}=&\frac{2BZX}{(ABC)^2},\\
\frac{dC}{dt}=&\frac{2CXY}{(ABC)^2},
\end{aligned}
\right .
\end{equation}
where
\begin{align*}
X=&(A-B)(3A+B), \\
Y=&(B-A)(3B+A), \\
Z=&-(A-B)^2.
\end{align*}

If $A_0=B_0$, then  $A=A_0$, $B=B_0$ and $C=C_0$. The geometry
stays flat all time. Without loss of generality, we assume that
$A_0>B_0$. Then  $A$ is increasing, while $B$ and $C$ are
decreasing.
\begin{lemma}Assume $A_0>B_0$. Then there is a finite time $T$
such that the solution exists on $[0,T)$ and
$$\lim_{t\rightarrow T} A(t)=\infty,\;\lim_{t\rightarrow T} B(t)=0
\mbox{ and }\lim_{t\rightarrow T} C(t)=0.$$
\end{lemma}
\begin{proof}
Notice that
$$X>\eta A^2,|Y|>\eta A^2 \mbox{ and }|Z|>\eta A^2$$ for some
constant $\eta\in (0,\infty)$. Hence $$-\ddt C>\frac{2\eta^2
A^4}{A^2B^2C}>\eta_1>0.$$ It follows that the solution exists up
to some finite time $T\in (0,\infty)$. Next
$$\ddt \ln
(AB)=\frac{2Z}{A^2B^2C^2}(X+Y)=\frac{-4(A-B)^4}{A^2B^2C^2} <0,$$
hence $AB$ is decreasing. Similarly,
$$\ddt \ln
(AC)=\frac{2Y}{A^2B^2C^2}(X+Z)<0,$$ so $AC$ is also decreasing.

We first show that $\lim_{T} C=0.$ Assume not. Then $\lim_{T}
C=C(T)\in (0,\infty)$. Since $AC$ is decreasing, we get $\lim_{T}
A=A(T)<\infty$ while $\lim_{T} B=0$. As $t\rightarrow T$, we have
$$X\sim 3A^2,~Y\sim -A^2, ~Z\sim -A^2,$$ hence $$\ddt
B^2=\frac{4ZX}{A^2C^2} \sim -\eta$$ for some positive fine $\eta$.
It follows that $B\sim \sqrt{\eta(T-t)}$ and
$$\ddt A=\frac{2AYZ}{(ABC)^2} \sim \frac{\eta}{T-t},$$
for a different $\eta\in (0,\infty)$. This shows that $A(t)\sim
-\eta \ln (T-t) \rightarrow \infty$ as $t\rightarrow T$, which is
a contradiction.

Next we show that $\lim_{T} A=\infty.$ Assume that we have
$\lim_{T} A=A(T)<\infty$. If further more we have $\lim_{T}
B=B(T)>0$ then, since $$\ddt C^2=\frac{4XY}{(AB)^2}\sim -\eta,$$
and thus $C(t)^2\sim \eta (T-t)$, just as above, we can show that
$A\rightarrow \infty$, this is a contradiction.

If we have $\lim_{T} B=0$, then $$B^2 \ddt C^2=\frac{4XY}{A^2}
\sim -12A(T)^2,$$ and
$$C^2 \ddt B^2=\frac{4XZ}{A^2} \sim -12A(T)^2.$$ Hence
$\ddt (B^2C^2)\sim -24A(T)^2$ and  thus $B^2C^2 \sim 24
A(T)^2(T-t)$. This again leads to $A\rightarrow \infty$, which is
a contradiction. This proves that $\lim_T A=\infty$.

Since $AB$ is decreasing, so we have proved that $\lim_{T}
A=\infty$ and  $\lim_{T} B=\lim_{T} C=0$.\end{proof}

\begin{lemma}Assume that $A_0>B_0$. Then there are constants
$\eta_1,\eta_2,\eta_3\in (0,\infty)$ such that
$$
A \sim  \eta_1 (T-t)^{-\frac{1}{14}},\;\;
B \sim  \eta_2 (T-t)^{\frac{3}{14}},\;\;
C \sim  \eta_3 (T-t)^{\frac{3}{14}}.
$$
\end{lemma}
\begin{proof}
As
$$\ddt (A^3B)=-\frac{16AZ}{C^2}(A-B)>0 \mbox{ and }
\ddt(A^3C)=\frac{8AY}{BC}(A-B)<0,$$ it follows that $A^3B$ is
increasing and $A^3C$ decreasing. We also have
$$\ddt (A^{2}B)=\frac{2A^{2}BZ}{A^2B^2C^2}(X+2
Y).$$ As $t\rightarrow T$, $X+2 Y =(A-B)(A-5B)>0.$ This implies
that $A^2B$ is bounded from above.  Using the same technique as in
the proof of Claim \ref{a3b}, we obtain that there exists $\eta\in
(0,\infty)$ such that
$\lim_{T} A^3B=\eta.$
Finally, observe that $$\ddt\ln(A^3C)\sim -8\frac{A}{BC^2} \mbox{
and } \ddt B\sim -6\frac{A^2}{BC^2}.$$ Since $B>0$ and
$\lim_TA=\infty$, it follows that $\int^T A/(BC^2)\le \int^T
A^2/(BC^2) <\infty$. Hence
$\lim_{T}
A^3C=\eta'\in (0,\infty).$ These asymptotic behaviors of $A^3B$
and $A^3C$ imply that
$$
\frac{dA}{dt}\sim  \eta A^{15},\;\;
\frac{dB}{dt}\sim  -\eta B^{-\frac{11}{3}},\;\;
\frac{dC}{dt}\sim -\eta C^{-\frac{11}{3}}.
$$
The desired asymptotic results follow.\end{proof}

\begin{theorem}\label{th-E22}
Let $g(t)$ be the solution of {\em (+XCF)} on $E(2)$ with given
initial data $A_0,B_0,C_0>0$ in a Milnor frame
 $f_1,f_2,f_3$ as above.\begin{itemize}
\item Assume that $A_0=B_0$. Then $g(t)=g_0=A_0f^1\otimes
f^1+A_0f^2\otimes f^2+C_0f^3\otimes f^3$. \item If $A_0>B_0$, then
there exists a finite time $T>0$, such that $g(t)$ exists on
$[0,T)$ and,  as $t\to T$,
$$
A \sim  \eta_1 (T-t)^{-\frac{1}{14}},~ B \sim \eta_2
(T-t)^{\frac{3}{14}},~and~ C \sim  \eta_3 (T-t)^{\frac{3}{14}}
$$ for some constants $\eta_1,\eta_2,\eta_3\in (0,\infty)$.
\end{itemize}
\end{theorem}

\section{The cross curvature flow on $\mbox{SL}(2,\mathbb R)$}

Given a left-invariant  metric $g_0$ on $\mbox{SL}(2,\mathbb R)$,
we fix a Milnor frame $\{f_i\}_1^3$ such that
$$[f_2,f_3]=-2f_1, \;\;[f_3,f_1]=2f_2,\;\;[f_1,f_2]=2f_3$$
and $$g_0=A_0f^1\otimes f^1+B_0f^2\otimes f^2+C_0f^3\otimes f^3.$$
The sectional curvatures are
\begin{align*}
K(f_2 \wedge f_3)&=\frac{1}{ABC}(-3A^2+B^2+C^2-2BC-2AC-2AB),\\
K(f_3 \wedge f_1)&=\frac{1}{ABC}(-3B^2+A^2+C^2+2BC+2AC-2AB),\\
K(f_1 \wedge f_2)&=\frac{1}{ABC}(-3C^2+A^2+B^2+2BC-2AC+2AB).
\end{align*}

Recall that the Lie algebra $\mbox{sl}(2,\mathbb R)$ of
$\mbox{SL}(2,\mathbb R)$ can be realized as the space of two by
two real matrices with trace $0$. A basis of this space is
$$W=\left(\begin{array}{cc}0&-1\\1&0\end{array}\right),\;\;
H=\left(\begin{array}{cc}1&0\\0&-1\end{array}\right),\;\;
V=\left(\begin{array}{cc}0&1\\1&0\end{array}\right).$$ These
satisfy $$[H,V]=-2W,\; [W,H]=2V,\; [V,W]=2H.$$ This means that
$(W,V,H)$ can be taken as a concrete representation of the above
Milnor basis $(f_1,f_2,f_3)$. In particular, $f_1$ corresponds to
rotation in $\mbox{SL}(2,\mathbb R)$. Note further that exchanging
$f_2,f_3$ and replacing $f_1$ by $-f_1$ produce another Milnor
basis. This explains the $B,C$ symmetry of the formulas above.

As for the other cases, the forward behavior of (XCF) was studied
in \cite[sect. 5]{cnsc1}. However, in the $\mbox{SL}(2,\mathbb R)$
case, the description of this forward asymptotic behavior  in
terms of convergence of metric spaces becomes quite intricate and
we will only make some simple remarks. The reader can consult
\cite{dg07} for a careful analysis using groupoid techniques.

 The backward behavior is studied below in the form of the
forward behavior of the positive cross curvature flow. Because of
the symmetry between $f_2$ and $f_3$ noted above, we can assume
that $B_0\ge C_0$.
\begin{theorem}\label{th-sl21}
Let $g(t)=A(t)f^1\otimes f^1+B(t)f^2\otimes  f^2  + C(t)f^3\otimes
f^3$, $t\in (-T_b,T_f)$, be a maximal solution of the negative
cross curvature flow on a complete locally homogeneous manifold of
type $\mbox{\em SL}(2,\mathbb R)$. Then $T_b<\infty$ whereas
$T_f=\infty$ if $B_0=C_0$ and $T_f<\infty$ otherwise. Assume that
$B_0\ge C_0$ and set
$$Q=\{(a,b,c)\in \mathbb R^3: a>0,
b\ge c>0\}$$ and
$$\bar{g}(t)=\frac{C_0}{C(t)} g(t).$$
There is a partition of $Q$ into subsets $S_0,Q_1,Q_2$ with
$Q_1,Q_2$ connected and, as $t$ tends to $-T_b$:
\begin{enumerate}
\item  If $(A_0,B_0,C_0)\in Q_1$, $(M,\overline{g}(t))$ converges
uniformly to the sub-Riemannian metric space $(M, bf_2\otimes f_2+
cf_3\otimes f_3)$ for some $b,c\in (0,\infty)$.

\item  If $(A_0,B_0,C_0)\in Q_2$, $(M,\overline{g}(t))$ converges
uniformly to the sub-Riemannian metric space $(M, af_1\otimes f_1+
cf_3\otimes f_3)$ for some $a,c\in (0,\infty)$.

\item If $(A(t),B(t),C(t))\in S_0$ for all $t\in (-T_b,0]$ then
$A(t)$ tends to $0$ whereas $B(t)$ and $C(t)$ converge towards the
same finite constant. If $M$ is compact, $(M,\overline{g}(t))$
converges in the Gromov-Hausdorff sense to a compact surface of
constant negative curvature.
\end{enumerate}
\end{theorem}
\begin{remark} In the forward direction, if $B_0=C_0$
and $M$ is compact then $(M,\overline{g}(t))$ converges in the
Gromov-Hausdorff sense to a compact  surface of constant negative
sectional curvature. See \cite{cnsc1} and \cite{dg07}.
\end{remark}
\begin{remark} The cases (1)-(2) of Theorem \ref{th-sl21} are
somewhat symmetric. Case (1) occurs when $A_0$ is large compared
to $B_0-C_0$. Case (2) occurs when $A_0$ is small compared to
$B_0-C_0$. Case (3) is of a completely different nature and it is
not even entirely clear, a priori, that it occurs at all. In a
forthcoming work \cite{cgsc08}, we will show that $Q_1\cup Q_2$ is
a dense open set and that $S_0$ is an hypersurface separating
$Q_1$ from $Q_2$. This however requires different techniques that
those used in this paper.
\end{remark}

The rest of this section is devoted to the proof of Theorem
\ref{th-sl21}. A much more precise statement is given in Theorem
\ref{th-sl2}. As in earlier sections, we focus on the forward
behavior of solutions of (+XCF). Using the sectional curvatures
given above, writing
$$g=  Af^1\otimes f^1+Bf^2\otimes f^2+Cf^3\otimes f^3$$
for the solution of the flow (+XCF) with initial data $g_0$,
$A,B,C$ satisfy the equations

\begin{equation}\label{pdesl2}
\left \{
\begin{aligned}
\frac{dA}{dt}=&\frac{2AF_2F_3}{(ABC)^2},\\
\frac{dB}{dt}=&\frac{2BF_3F_1}{(ABC)^2},\\
\frac{dC}{dt}=&\frac{2CF_1F_2}{(ABC)^2},
\end{aligned}
\right .
\end{equation}
where
\begin{align*}
F_1=&-3A^2+B^2+C^2-2BC-2AC-2AB, \\
F_2=&-3B^2+A^2+C^2+2BC+2AC-2AB, \\
F_3=&-3C^2+A^2+B^2+2BC-2AC+2AB.
\end{align*}

Without loss of generality we may assume that $B_0 \geq C_0$. Then
$B \geq C$ as long as a solution exists and
$$
F_3=(B-C)(2A+B+3C)+A^2 \geq A^2>0.
$$
Observe also that
$F_1+F_2<0$ so at least one of the quantities $F_1,F_2$ is negative.

Let $a=A/B$ and $c=C/B$.
\begin{lemma}\label{lem-f2n}
Suppose that  $a_0=A_0/B_0$ and $c_0=C_0/B_0$ satisfy
\begin{equation}\label{f2n}
a>1-c+2\sqrt{1-c}
\end{equation}
Then $a$ and $c$ satisfy {\em (\ref{f2n})} as long as a solution
exists. Moreover, in this case, {\em (\ref{f2n})} is equivalent to
$F_2>0$
 and thus implies that  $F_1<0$.
\end{lemma}

\begin{proof}
As in \cite[Sect. 5, Lemma 2]{cnsc1}, we have
$$F_2=(A-(B-C+2\sqrt{(B-C)B}))(A-(B-C-2\sqrt{(B-C)B})).$$
Since $ B-C-2\sqrt{(B-C)B}\le 0$, it follows that
 (\ref{f2n}) is equivalent to $F_2>0$. Observe that
$$
\left.\frac{dA}{dt}\right|_{F_2=0}=0,\quad
\left.\frac{dB}{dt}\right|_{F_2=0}<0\quad\mbox{ and }
\left.\frac{dC}{dt}\right|_{F_2=0}=0.
$$
It follows that
\begin{equation}\label{df2dt}
\left.\frac{dF_2}{dt}\right|_{F_2=0}>0.
\end{equation}
To prove that $F_2(t)>0$ we argue by contradiction. Suppose $t_0$
is the first time such that $F_2(t_0)=0$. Since $F_2(0)>0$, we
know that $\ddt F_2(t_0)\leq 0$, which contradicts (\ref{df2dt}).
Therefore $F_2(t)>0$, which is equivalent to (\ref{f2n}). This
completes the proof of the lemma.
\end{proof}

The next lemma is very similar to the previous one and we omit the proof.
\begin{lemma}\label{lem-f1n}
Suppose that $a_0=A_0/B_0$ and $c_0=C_0/B_0$ satisfy
\begin{equation}\label{f1n}
a<\frac13(2\sqrt{1-c+c^2}-1-c).
\end{equation}
Then $a$ and $c$ satisfy and {\em (\ref{f1n})} as long as a
solution exists. Moreover, {\em (\ref{f1n})} is equivalent to
$F_1>0$ and thus implies $F_2<0$.
\end{lemma}
Observe that, if at any time $t_0$, the solution satisfies either
(\ref{f2n}) or (\ref{f1n}) then that inequality will be satisfied
at all later time. We will consider three cases: The case
where (\ref{f2n}) is satisfied (at time $t=0$ or, in fact at a later time),
the case where (\ref{f1n}) is satisfied (at time $t=0$ or, in fact at
a later time), and the remaining case where neither (\ref{f2n}) nor
(\ref{f1n})
is satisfied as long as the solution exists.

{\bf Case 1: Inequality (\ref{f2n}) is satisfied.}  Recall that
this is equivalent to say that  $F_2>0$. Moreover, we must have $F_1<0$.

\begin{lemma} \label{f13a2} Assume that $B\ge C$ and {\em
(\ref{f2n})}
holds. Then $A$ is increasing, $B$ and $C$ are decreasing and
$A+C\ge B$. Furthermore $a$,  $c$ and $F_2/B^2$
are non-decreasing and  $|F_1|>3A^2$.
\end{lemma}

\begin{proof}  The monotonicity of $A$, $B$ and $C$ follows
from the fact that $F_2>0$ which also easily imply $A+C\ge B$.
Further
\begin{equation}\label{dlna}
\frac{d\ln(A/B)}{dt}=\frac{8(A+B)}{(ABC)^2}F_3(A+C-B),
\end{equation}
and
\begin{equation}\label{dlnc}
\frac{d\ln(C/B)}{dt}=\frac{8(B-C)}{(ABC)^2}(A+B+C)|F_1|,
\end{equation}
are non-negative. As
$$\frac{F_2}{B^2}=(\frac{A+C}{B}-1)^2+4\frac{C}{B}-4,$$
and
$$|F_1|= 3A^2+2(A+C)B+2AC-B^2-C^2.$$
The lemma follows.
\end{proof}

\begin{lemma}
Assume that $B\ge C$ and {\em (\ref{f2n})} is satisfied. Then
there exists a finite time $T$ such that
$$\lim _{t\rightarrow T}A(t)=\infty,\;\;
 \lim _{t\rightarrow T}B(t)=\lim _{t\rightarrow T}C(t)=
 0.$$
 Moreover, there exists $\eta_0\in (0,\infty)$ such that
 $\lim_TB/C=\eta_0$.
\end{lemma}
\begin{proof}
By Lemma \ref{f13a2}, there exists $\eta\in (0,\infty)$ such that
\begin{equation}
\ddt C^2=-4 \frac{|F_1|}{A^2} \times \frac{F_2}{B^2}<-12 \eta .
\end{equation}
Hence there exists a finite time $T$ such that the flow
exists only up to $T$ and  either  $\lim_TC=0$ or
$\lim_{T} A=\infty$.

Observe that $C/B\le 1$ is non-decreasing.
Hence  $\lim_T B/C\in (0,\infty)$.
We first show that $\lim_TC=0$, which implies that $\lim_TB=0$.
Otherwise we have $\lim_{T} A=\infty$ and $\lim_T
C=C(T) >0$. Note that $\lim_TB=B(T)\ge C(T)>0$.
Since $$\ddt A=\frac{2}{A} \times \frac{F_2}{B^2} \times
\frac{F_3}{C^2},$$ and
$F_2 \sim F_3 \sim A^2$ as $t\rightarrow T$,
there exists a constant $\eta >0$, such
that $$\ddt (A^{-2}) \sim -\frac{1}{\eta}.$$ Hence $$A^2 \sim
\frac{\eta}{T-t}.$$  As $|F_1|\sim 3A^2$,
$$\ddt C^2=-4 \frac{|F_1|}{A^2}
\frac{F_2}{B^2} \sim -\frac{12\eta}{B(T)^2}\times\frac{1}{T-t}.$$ This
contradicts the fact that $C$ is decreasing with $C(T)>0$.

Next, we prove that $\lim_TA=\infty$.
Observe that, for any $\alpha>0$, we have
\begin{equation}\label{alphaAB}
\ddt (A^{\alpha}B)=\frac{2A^{\alpha}BF_3}{(ABC)^2}(\alpha
F_2+F_1),
\end{equation} and
\begin{equation}\label{alphaAB'}
\alpha F_2+F_1=(\alpha
-3)A^2+(1-3\alpha)B^2+(1+\alpha)C^2+2(\alpha-1)(B+A)C-2(1+\alpha)AB.
\end{equation}
If $\alpha >3$, then $\lim_{T} (\alpha
F_2+F_1)>0$ and thus
$A^{\alpha}B$ increasing. As $\lim_TB= 0$, it follows that
$\lim_{T} A=\infty$, as desired. \end{proof}

\begin{lemma} Assume that $B\ge C$ and {\em (\ref{f2n})} holds. There are
constants $\eta_1,\eta_2,\eta_3\in (0,\infty)$ such that,
as $t$ tends to $T$,
$$A(t) \sim \eta_1 (T-t)^{-1/14},\;\;B(t) \sim \eta_2 (T-t)^{3/14}
\mbox{ and }\;C(t) \sim \eta_3 (T-t)^{3/14}.$$
\end{lemma}
 \begin{proof}
For $\alpha =3$, $ (\alpha F_2+F_1)=-4[(B-C)(A+2B+C)+AB]<0.$
By (\ref{alphaAB}) it follows that $A^3B$ is decreasing.
Further, $$\ddt
\ln(A^3B)=-8\frac{F_3}{(ABC)^2}[(B-C)(A+2B+C)+AB] \sim
-8\frac{A}{B^2C^2}(2B-C)$$
and
\begin{equation}\label{dlncasy}
\ddt B \sim -6\frac{A^2}{BC^2} \mbox{ as } t\rightarrow T .
\end{equation}
As  $B>0$, we must have
$$\int^{T} \frac{A^2}{BC^2} < \infty.$$
It follows that $$\int^{T} \frac{A}{B^2C^2} (2B-C)
< 2\int^{T} \frac{A}{BC^2} < \int^{T} \frac{A^2}{BC^2} <
\infty.$$ Hence there exists $\eta\in (0,\infty)$ such that
$\lim_{T} A^{3}B=\eta.$
Now, as $t$ tends $T$,
$$\ddt A \sim \eta \frac{A^5}{A^2B^2C^2} \sim \eta A^{15}.$$
The asymptotic for $A$ follows as well as those for $B$ and $C$.
\end{proof}

{\bf Case 2: Inequality (\ref{f1n}) is satisfied.} This is equivalent
to $F_1>0$ and  implies that $F_2<0$.
Since $B\ge C$ and $F_1>0$, we have $A+C<B$. Since $F_3>0$, we
have $A$ is decreasing, $B$ is increasing and $C$ is decreasing,
hence both $a$ and $c$ are decreasing. Further,
$$\frac{F_1}{B^2}=1+\frac{C^2}{B^2}-3\frac{A^2}{B^2}-2\frac{C}{B}-2\frac{A}{B}
\frac{C}{B}-2\frac{A}{B}=(1-c)^2-3a^2-2ac-2a.$$ Since both $a>0$
and $c>0$ are decreasing, $\frac{F_1}{B^2}$ is
increasing, hence as long as the solution exists, we have
$$\frac{F_1}{B^2} \geq  \frac{F_1(0)}{B_0^2}>0.$$
We also have
$$|F_2|=3B^2-A^2-C^2-2BC-2AC+2AB\geq 2B^2-2BC+2AB>4A^2.$$
Hence there exists $\eta\in (0,\infty)$ such that
$$-\ddt
C^2=4\frac{F_1}{B^2} \frac{F_2}{A^2} >\eta.$$ This implies that
the solution can only exist up to some  finite time $T$ at which
at least one of the following must happen: $\lim_T A=0$,  $\lim_T
C=0$ or $\lim_T B=\infty$.
\begin{lemma}Assume that $B\ge C$ and {\em (\ref{f1n})} holds. Then
there exists $T\in (0,\infty)$ such that the solution exists on $[0,T)$
and
$$\lim_{t\rightarrow T}A(t)=\lim_{t\rightarrow T}C(t)=0 \mbox{ and
} \lim_{t\rightarrow T}B(t)=\infty.$$
Further, there exists $\eta_0\in (0,\infty)$ such that $\lim_T C/A=\eta_0$.
\end{lemma}
\begin{proof}As
\begin{equation}\label{f1nAC}
\ddt \ln
(C/A)=8 \frac{F_2}{A^2B^2C^2} (C+A)(C-A-B)>0,\end{equation}
$C/A$ is increasing. Hence, if $\lim C=0$,
then we must have $\lim A=0$.

We first show that
 $\lim_{t\rightarrow T} \frac{B}{C} =\infty$, i.e.,
 $\lim_{t\rightarrow T} c=0$. Indeed,
if $\lim_T C=0$, since $B$ is increasing, and the desired result follows.
If $\lim_T C=C(T)>0$, it suffices to show that $\lim B=\infty$.
Assume  instead that $\lim_T B =B(T)<\infty$. Then we must have
$\lim_TA=0$ and, as
$t\rightarrow T$,
$$F_2 \sim -3B^2+C^2+2BC \mbox{ and } F_3 \sim
-3C^2+B^2+2BC.$$ Hence there exists $\eta\in (0,\infty)$ such that
$$\ddt A^2 \sim -4\frac{(B-C)^2 (B+3C)
(3B+C)}{B^2C^2}\sim \eta.$$ It follows that $A^2 \sim
\eta (T-t)$ as $t$ tends to $T$.  Also, as $t$ tends to $T$,
$$F_1 \sim (B-C)^2,$$ so $$\ddt C \sim -2 \frac{1}{\eta
(T-t)} \frac{(B-C)^3(3B+C)}{B^2 C^2} \sim - \frac{\eta'}{T-t},$$
for some $\eta'\in (0,\infty)$. This contradicts  the fact that
$C>0$. So $\lim_T B=\infty$ and $\lim_T B/C=\infty$.

Next we show that
$\lim_T C=0$ ($\lim_T A=0$ follows).
If not, we have $\lim_T C=C(T)>0$ and this implies that $\lim_T B=\infty$.
Then, since $A$ is decreasing, we have
$$F_1 \sim B^2,\;\; F_2 \sim -3B^2 \;\mbox{ and }\;F_3 \sim B^2.$$
So
$$\ddt \ln B \sim \frac{2}{C(T)^2}\frac{B^2 }{A^2}
\mbox{ and } \ddt C \sim -\frac{6}{C(T)}\frac{B^2 }{A^2}.$$
Since $\lim_T B =\infty$, we must have $\int^T \frac{B^2}{A^2}
=\infty.$ This contradicts the fact that  $C>0$. So
we have $\lim_T C=0$ as desired.

By (\ref{alphaAB})-(\ref{alphaAB'}), if $\alpha<\frac13$ and  $t$ is
close enough to $T$, $A^{\alpha}B$ is increasing, but $\lim A=0$,
so $\lim B=\infty$.

Finally, we prove    that $\lim_T \frac{C}{A} \in (0,\infty)$.
We already know from (\ref{f1nAC})
that  $C/A$ is increasing.
On the one hand,
we have $$\ddt \ln \frac{C}{A} \sim \frac{24B}{A^2C^2}
(A+C)\le \eta \frac{B^2}{A^2C}$$
for some $\eta\in (0,\infty)$. On the other hand, $$-\ddt
C \sim 6 \frac{B^2}{A^2C}.$$ Hence $$\int^T \frac{B^2}{A^2C}
<\infty.$$ It follows that $\ln \frac{C}{A}$ is bounded from above.
\end{proof}

\begin{lemma} Assume that $B\ge C$ and {\em (\ref{f1n})} holds.
There are constants $\eta_1,\eta_2,\eta_3\in (0,\infty)$ such that,
as $t$ tends to $T$,
$$A(t) \sim \eta_1 (T-t)^{3/14},\;\;B(t) \sim \eta_2 (T-t)^{-1/14}
\mbox{ and }\;C(t) \sim \eta_3 (T-t)^{3/14}.$$
\end{lemma}
\begin{proof}
We first show that $\lim_T A^{\frac13}B\in(0,\infty)$.
By (\ref{alphaAB})-(\ref{alphaAB'}) with  $\alpha=\frac13$, we have
$$\ddt
(A^{\frac13}B)=\frac{2A^{\frac13}BF_3}{(ABC)^2}(\frac13
F_2+F_1)\sim -\frac83\frac{A^{\frac13}BF_3}{A^2B^2C^2}(BC+2AB
)<0,$$ for $t$ close enough to $T$, so we only need to show that
$\lim_T A^{\frac13}B>0$.
We have
$$\ddt
\ln (A^{\frac13}B)\sim -\frac83\frac{F_3}{A^2B^2C^2}(BC+2AB)\sim
-\eta \frac{B}{A^2C},$$
for some constant $\eta\in (0,\infty)$ (here, we use that
$\lim_T C/A\in (0,\infty)$). Further,
$$\ddt A=\frac{2AF_2F_3}{(ABC)^2}\sim -6\frac{B^2}{AC^2}.$$ Since
$A$ is bounded, we must have  $\int^T \frac{B^2}{AC^2}<\infty$. Hence
$\int^T \frac{B}{A^2C}<\infty$ as well. This implies
$\lim_T A^{\frac13}B=\eta>0$ as desired.
It follows that $AB^3$ and $CB^3$ have positive finite limits as $t$ tends
to $T$ and we can proceed as
in case 1 to obtain the asymptotic of $A$, $B$ and
$C$ when $t$ tends to $T$.\end{proof}

{\bf Case 3: Neither (\ref{f2n})
nor  (\ref{f1n}) are ever satisfied along the flow}
Assume the solution exists on the interval $[0,T)$
($T$ can be $\infty$ here).
The third case is the case when
\begin{equation}\label{case3}\forall\,t\in [0,T),\;\;
 \frac{1}{3}\left(2\sqrt{1-c+c^2}-1-c\right)\le a\le 1-c+2\sqrt{1-c}.
 \end{equation}
Recall that $(\ref{f2n})$ is equivalent to $F_2>0$  and (\ref{f1n})
is equivalent to $F_1>0$. Hence we have $F_1,F_2\le 0$.
Since $F_1+F_2<0$, at least one of them is strictly negative
in this case.  We first notice
that $F_3\geq A^2>0$. It follows that
$A$ and $B$ are non-increasing, and $C$ is non-decreasing. So we have
$$C_0\leq C \leq B \leq B_0$$
and $B,C$ have finite positive limits when $t$ tends to $T$.
\begin{lemma} Assume  $B\ge C$.
\begin{itemize}\item If $a_0=A_0/B_0,c_0=C_0/B_0$ satisfy
\begin{equation}\label{f3n} a\ge 1-c.
\end{equation}
then this inequality is satisfied for all $t\in [0,T)$ and there
exists a time $t_1\in [0,T)$ such that  {\em (\ref{f2n})} holds
for all $t\in (t_1,T)$. \item If {\em (\ref{case3})} holds then,
for all $t\in [0,T)$, $A+C<B$, that is, $a< 1-c$.
\end{itemize}
\end{lemma}
\begin{proof}
Observe that for any $c\in [0,1]$,
$$1-c \ge \frac{(1-c)^2}{2\sqrt{1-c+c^2}+1+c}=
\frac13(2\sqrt{1-c+c^2}-1-c).$$ By Lemma \ref{lem-f1n}, this shows
that $a\ge 1-c$ implies $F_1\le 0$. Assume that $A+C\geq B$ at
some time $t$. Then, at that time, (\ref{dlna}) and (\ref{dlnc}) show that
$\ddt \ln (A/B) \geq 0$ and  $\ddt \ln (C/B)\geq 0$. This proves that the
inequality $a\ge 1-c$ is preserved by the flow.

Assume now that (\ref{case3}) holds. Then both $B$ and $C$ are
monotone and have positive finite limits as $t$ tends to $T$. If
there exists a time $t_0$ such that $a\ge 1-c$ at time $t_0$, then
$a=A/B$ is non-decreasing for $t\ge t_0$. This means that the
solution must exist for all time, i.e., $T=\infty$. However,
$$|F_1|=3A^2-(B-C)^2+2AC+2AB\geq 2A(A+B+C),$$
so $\ddt B<-\eta$ with $\eta$ a positive constant. This
contradicts $T=\infty$ and thus shows that $A+C<B$, for all $t\in
[0,T)$.

Finally, since $a\ge 1-c$ is preserved, implies $F_1\ge 0$, and is
incompatible with (\ref{f3n}), it follows that if $a_0\ge 1-c_0$
then there exists $t_1$ such that $F_2>0$, that is (\ref{f2n})
holds (it then holds for all $t>t_1$ by Lemma \ref{lem-f2n}).
\end{proof}

\begin{lemma}Assume  $B\ge C$ and {\em (\ref{case3})}. Then we have
$$\lim_{t\rightarrow T}A(t)=0\;\mbox{ and }\;\;
\lim_{t\rightarrow T}B(t)=\lim_{t\rightarrow T}C(t)=k\in
(0,\infty).$$
\end{lemma}
\begin{proof} Observe that $F_2+F_3=2(A+B-C)(A-B+C)<0$. It follows that
$|F_2|>F_3\geq A^2$. If $\lim_T A=A(T)>0$ then we must have
$T=\infty$. However, $\ddt A<-\eta$ with $\eta$ a positive
constant so the solution can only exist up to some finite time, a
contradiction. This shows that $\lim_TA=0$. Since $0\le B-C\le A
$, we also have  $\lim_T(B-C)=0$.
\end{proof}

\begin{lemma}Assume  $B\ge C$ and {\em (\ref{case3})}. Then $T< \infty$,
and  $$\lim_{t\rightarrow T}\frac{2kA}{(B-C)^2}=1.$$
\end{lemma}
\begin{proof}
First we observe that the inequality $4AB\ge (B-C)^2$ follows
easily from the fact that
$$a\geq \frac13(2\sqrt{1-c+c^2}-1-c)=
\frac{(1-c)^2}{2\sqrt{1-c+c^2}+1+c}.$$

For any positive number $\alpha$, we have
$$\ddt \ln A=\frac{2}{A^2B^2C^2} F_2 F_3 \mbox{ and }
\ddt \ln (B-C) =\frac{2}{A^2B^2C^2} F_1 Y.$$ It follows that
\begin{equation}\label{A/BC}
\ddt \ln \frac{A}{(B-C)^{\alpha}}=\frac{2}{A^2B^2C^2}
(F_2F_3-\alpha F_1Y),\end{equation}
 where
$Y= A^2+B^2+C^2+6BC+2AB+2AC$. Since $B-C>A$ and $\lim_T A=0$, as
$t$ tends to $T$, we have
$$F_2=-(2A+C+3B)(B-C)+A^2 = -4k (B-C) +O((B-C)^2),$$
$$F_3=(B-C)(B+3C+2A)+A^2 = 4k(B-C)+O((B-C)^2),$$
$$-F_1=2AB+2AC-(B-C)^2+3A^2=4kA-(B-C)^2+O(A(B-C)),$$
and $$Y=8k^2+O((B-C)).$$ In this estimates, we have used the fact
that $0\le k-C,\,B-k\le B-C$. This yields
\begin{equation}\label{FFFY}
F_2F_3-\alpha F_1Y = 8\alpha k^2[4k A-(1+(2/\alpha))(B-C)^2]
+O[A(B-C)+(B-C)^3].\end{equation}
 We first show that for any $\epsilon>0$ there exists $t_0$ such
that
$$\forall\,t\in [t_0,T),\;\; \frac{4kA}{(B-C)^2} \ge 2-\epsilon.$$
Otherwise, there exists $\epsilon>0$ such that for any $t_0$ there
exists  $t_1\in (t_0,T)$  such that $\left.4kA/(B-C)^2\right|_{t=
t_1} <2-\epsilon$. Taking $t_0$ large enough, $\alpha=2$ in
(\ref{A/BC}) and using (\ref{FFFY}), it follows that we must have
 $4kA/(B-C)^2 < 2-\epsilon$ for all $t\in [t_1, T)$. Now, taking
$\alpha=2+\delta(\epsilon)$ in (\ref{A/BC}), we obtain that
$A/(B-C)^{2+\delta(\epsilon)}$ is decreasing if $2-\epsilon <
1+\frac{2}{2+\delta(\epsilon)}$. This contradicts the fact that
$A/(B-C)^2 \ge 1/4B$.

Now, as $t$ tends to $T$,
$$|F_1|=2AC+2AB-(B-C)^2+3A^2\ge kA.$$
Hence, there exists $\eta\in (0,\infty)$ such that
$$\ddt C =\frac{2}{B^2C^2} \frac{|F_1|}{A} \times \frac{|F_2|}{A}
>\eta.$$
As $C$ is bounded from above by $B_0$, this implies $T<\infty$.

Next we show that there exists $t_0\in ([0,T)$ such that
$$\forall \,t\in [t_0,T),\;\;\frac{kA}{(B-C)^2} < 1.$$
Otherwise, for any $t_0$ there exists  $t_1\in
(t_0,T)$ such that $\left.kA/(B-C)^2\right|_{t=t_1}\ge 1$. Using
(\ref{A/BC})-(\ref{FFFY}), it is easy to see the $kA/(B-C)^2 \ge 1
$ for all $t\in [t_1,T)$. Now take $\alpha=1$ in
(\ref{A/BC})-(\ref{FFFY}) to see that  $A/(B-C)$ is increasing on
$(t_2,T)$ for some $t_2\in (t_1,T)$. It follows that there exists
$\eta\in (0,\infty)$ such that
 $$-\ddt \ln
A\sim \eta \frac{(B-C)^2}{A^2} <\eta.$$ As $T<\infty$, this
contradicts $\lim_{T} A=0$.

We now show that for any $\epsilon>0$ there exists $t_0$ such that
$$\forall\,t\in [t_0,T),\;\; \frac{4kA}{(B-C)^2} \le 2+\epsilon.$$
Otherwise, there exists $\epsilon>0$ such that for any $t_0$ there
exists  $t_1\in (t_0,T)$  such that $\left.4kA/(B-C)^2\right|_{t=
t_1} >2+\epsilon$. Taking $t_0$ large enough, $\alpha=2$ in
(\ref{A/BC}) and using (\ref{FFFY}), it follows that we must have
 $4kA/(B-C)^2 > 2+\epsilon$ for all $t\in [t_1, T)$.
Now using $\alpha=2-\delta(\epsilon)$ in
(\ref{A/BC})-(\ref{FFFY}), we find that
$A/(B-C)^{2-\delta(\epsilon)}$ is non-decreasing if $2+\epsilon >
1+\frac{2}{2-\delta(\epsilon)}$. This contradicts $A/(B-C)^2 \le
3/4k$.\end{proof}

Using $2kA\sim (B-C)^2$ in the original ODEs (\ref{pdesl2}), we
easily find the asymptotic behaviors of $A,B$ and $C$. Namely, as
$t\rightarrow T$, we have
$$A\sim \frac{64}{k} (T-t),\;\;B\sim k+4\sqrt{2} \sqrt{T-t}
\mbox{ and } C\sim k-4\sqrt{2} \sqrt{T-t}.$$

Putting together the different cases yields the following
statement.
\begin{theorem}\label{th-sl2}
Let $g(t)$ be the solution of {\em (+XCF)} on $\mbox{\em
SL}(2,\mathbb R)$ with given initial data $A_0,B_0,C_0>0$, $B_0\ge
C_0$, in a Milnor frame $f_1,f_2,f_3$. Then there exists a finite
positive time $T_b$ such that $g(t)$ exists for all $t\in
[0,T_b)$. Moreover, we have:
\begin{enumerate}
\item If there is a time $t$ such that $A\ge B-C$ then there are
constant $\eta_i\in (0,\infty)$ such that, as $t\to T_b$,
$$ A \sim \eta_1 (T_b-t)^{-\frac{1}{14}},~ B \sim  \eta_2
(T_b-t)^{\frac{3}{14}},~and~ C \sim  \eta_3
(T_b-t)^{\frac{3}{14}}.
$$
\item If there is a time $t$ such that $A<\frac13
(2\sqrt{B^2-BC+C^2}-B-C)$ then there are constant $\eta_i\in
(0,\infty)$ such that, as $t\to T_b$,
$$ A \sim \eta_1 (T_b-t)^{\frac{3}{14}},~ B \sim \eta_2
(T_b-t)^{-\frac{1}{14}},~and~ C \sim  \eta_3
(T_b-t)^{\frac{3}{14}}.
$$
\item If for all $t$, $\frac13 (2\sqrt{B^2-BC+C^2}-B-C)\leq A <
B-C$, then as $t\to T_b$,
$$
A\sim \frac{64}{C(T_b)} (T_b-t),~ B\sim
C(T_b)+4\sqrt{2(T_b-t)},~and~ C\sim C(T_b)-4\sqrt{2(T_b-t)}
$$
with $C(T_b)\in (0,\infty)$.
\end{enumerate}
\end{theorem}
\begin{remark} A priori,
the statement above does not imply that the case when
$$\frac13 (2\sqrt{B^2-BC+C^2}-B-C)\leq A < B-C$$ holds for all
$t\in [0, T)$ really occurs since it is possible that no solutions
satisfy $\frac13 (2\sqrt{B^2-BC+C^2}-B-C)\leq A < B-C$ for all
$t\in [0,T_b)$. In \cite{cgsc08}, we will show by different
techniques that this case really occurs but only for initial
condition on an hypersurface. See the discussion below.
\end{remark}
\begin{remark}
In case (3), we have
$$\lim_{t \rightarrow T_b} K(f_2\wedge f_3)=-\frac{2}{C(T_b)}$$
and the two other sectional curvatures vanish. If $M$ is compact,
then there is collapse in the $f_1$ direction.
\end{remark}

It is quite clear that Theorem \ref{th-sl2} implies Theorem
\ref{th-sl21} and give a more precise and technical description of
the asymptotic behavior of (+XCF) on $\mbox{SL}(2,\mathbb R)$.
More precisely, the initial condition space
$$Q=\{(a,b,c)\in\mathbb R^3:
a>0, b\ge c>0\}$$ can be partitioned into $Q=Q_1\cup Q_2\cup S_0$
with
$$\{(a,b,c)\in Q: a\ge b-c\}\subset Q_1,$$
$$\{(a,b,c)\in Q: 3 a< 2\sqrt{b^2-bc +c^2} -(b+c)\}\subset Q_2, $$
and
$$S_0\subset  \{(a,b,c)\in Q:
( 2\sqrt{b^2-bc +c^2} -(b+c)) \le 3a< 3(b-c)\}.$$ By Theorem
\ref{th-sl2}(1)-(2), for initial condition in $Q_1$ (resp. $Q_2$)
the distance function associated to the metric $\overline{g}=
(C_0/C(t))g(t)$ clearly converges uniformly on compact sets to the
distance function associated with a sub-Riemannian structure of
the form $\gamma_2f_2\otimes f_2+ \gamma_3 f_3\otimes f_3$ (resp.
$\gamma_1 f_1\otimes f_1+\gamma_3 f_3\otimes f_3$). In fact, as
stated in  Theorem \ref{th-sl2}, it suffices that $g(t)$ enters
the region $Q_1$ or $Q_2$ at some time $t\ge 0$. Hence, we can
define $Q_1$ (resp. $Q_2$) to be the set of initial conditions in
$Q$ such that $g(t)$ enters $\{a\ge b-c\}$ (resp. $ \{3 a<
2\sqrt{b^2-bc +c^2} -(b+c)\}$) and $S_0=Q\setminus (Q_1\cup Q_2)$.
From this discussion it  seems very plausible that $S_0$ is simply
a hypersurface separating $Q_1$ and $Q_2$. However, we have not
been able to prove this by arguments similar to those used above.
In \cite{cgsc08}, we prove that $S_0$ is indeed a surface
separating the open sets $Q_1,Q_2$ by reducing the ODEs to a two
dimensional one and using the fact that the orbit structure of
$2$-dimensional ODEs can be understood much better than in higher
dimension.

\bibliographystyle{alpha}
\bibliography{bio}
\end{document}